\RequirePackage{fix-cm}
\documentclass[twocolumn]{svjour3}          
\smartqed  
\usepackage{dblfloatfix}
\usepackage{float}
\usepackage{graphicx}
\usepackage{amsmath}
\usepackage{lineno,hyperref}
\usepackage{times}
\usepackage{epsfig}
\usepackage{bm}
\usepackage{breakurl}
\usepackage{makecell}
\usepackage{graphicx, subfigure}
\usepackage{amssymb, amsfonts, euscript, mathrsfs}
\usepackage{subeqnarray}
\usepackage{cases}
\usepackage[sectionbib,round]{natbib}

\usepackage{overpic}
\usepackage{listings}
\usepackage{array}

\usepackage{fancybox}
\usepackage{picins}
\usepackage{picinpar} 
\usepackage[table]{xcolor}
\usepackage{wrapfig}
\definecolor{lightgray}{gray}{0.9}
\definecolor{lightblue}{rgb}{0.98,0.98,1.0}
\usepackage{subfigure}
\usepackage{multirow}
\usepackage{booktabs}
\usepackage{ulem}
\usepackage{comment}
\newcommand{\executeiffilenewer}[3]{%
\ifnum\pdfstrcmp{\pdffilemoddate{#1}}%
{\pdffilemoddate{#2}}>0%
{\immediate\write18{#3}}\fi%
}

\usepackage[ruled]{algorithm2e} 

\usepackage{multirow} 
\usepackage{xcolor}

\newcommand{\T}{\mathsf{T}}

\newcommand{\tr}{\textcolor[rgb]{1.0,0.0,0.0}}
\newcommand{\tb}{\textcolor[rgb]{0.0,0.0,1.0}}

\newcommand{\fuxm}[1]{\tb{[fuxm:#1]}}

\makeatletter
\DeclareRobustCommand\bigop[1]{%
  \mathop{\vphantom{\sum}\mathpalette\bigop@{#1}}\slimits@
}
\newcommand{\bigop@}[2]{%
  \vcenter{%
    \sbox\z@{$#1\sum$}%
    \hbox{\resizebox{\ifx#1\displaystyle.9\fi\dimexpr\ht\z@+\dp\z@}{!}{$\m@th#2$}}%
  }%
}
\makeatother

\newcommand{\stkout}[1]{\ifmmode\text{\sout{\ensuremath{#1}}}\else\sout{#1}\fi}

\graphicspath{ {./images/} }

\modulolinenumbers[5]
\journalname{Structural and Multidisciplinary Optimization}
\begin{document}\sloppy

\definecolor{mygreen}{rgb}{0,0.6,0}
\definecolor{mygray}{rgb}{0.5,0.5,0.5}
\definecolor{mymauve}{rgb}{0.58,0,0.82}
\lstset{ %
backgroundcolor=\color{lightgray},   
basicstyle=\footnotesize\ttfamily,   
columns=fullflexible,
breaklines=true,                 
captionpos=b,                    
tabsize=4,
commentstyle=\color{mygreen},    
escapeinside={\%*}{*)},          
keywordstyle=\color{blue},       
stringstyle=\color{mymauve}\ttfamily,     
frame=single,
rulesepcolor=\color{red!20!green!20!blue!20},
language=c++,
}

\title{An Optimized, Easy-to-use, Open-source GPU Solver for Large-scale Inverse Homogenization Problems}


\subtitle{}


\author{Di Zhang$^{1}$  \and
        Xiaoya Zhai$^{1*}$ \and
        Ligang Liu$^{1}$  \and
        Xiao-Ming Fu$^{1}$
}


\institute{
        \at
        $^1$School of Mathematical Sciences, University of Science and Technology of China, Hefei, China \\
        \at
        $^*$Corresponding Author: \email{xiaoya93@mail.ustc.edu.cn} (Xiaoya Zhai)
}

\date{Received: date / Accepted: date}

\maketitle

\begin{abstract}
We propose a high-performance GPU solver for inverse homogenization problems to design high-resolution 3D microstructures.
Central to our solver is a favorable combination of data structures and algorithms, making full use of the parallel computation power of today's GPUs through a software-level design space exploration.
This solver is demonstrated to optimize homogenized stiffness tensors, such as bulk modulus, shear modulus, and Poisson’s ratio, under the constraint of bounded material volume.
\tb{Practical high-resolution examples with $512^3 \approx 134.2$ million finite elements run in less than 40 seconds per iteration with a peak \tb{GPU} memory of 9 GB} \tb{on an NVIDIA GeForce GTX 1080Ti GPU.}
Besides, our GPU implementation is equipped with an easy-to-use framework with less than 20 lines of code to support various objective functions defined by the homogenized stiffness tensors.
Our open-source high-performance implementation is publicly accessible at \url{https://github.com/lavenklau/homo3d}.
\end{abstract}


\section{Introduction} \label{sec:Intro}
Microstructure design is fundamental in various applications, such as aerospace and biomedicine.
Topology optimization for inverse homogenization problems (IHPs)~\citep{sigmund1994materials} is a powerful and effective method to find optimal microstructures.
Many methods have been developed for the microstructure topology optimization, such as density-based method~\citep{aage2015topology,groen2018homogenization}, isogemetric topology optimization~\citep{gao2019topology,gao2020isogeometric}, bidirectional evolutionary structural optimization~\citep{huang2012evolutionary,huang2011topological}, and level set method~\citep{vogiatzis2017topology,li2018topology}.


We focus on high-resolution periodic 3D microstructure design via the density-based method.
High-resolution microstructures expand the search space of multi-scale structures, making their mechanical properties more likely to approach the optimal solution. 
In addition, a microstructure is periodically arranged in a macro-scale domain behaving like a material under the premise that the length scale of the microstructure is much smaller than that of the macrostructure based on the homogenization theory~\citep{bendsoe1988generating,suzuki1991homogenization,nishiwaki1998topology}.

We aim to use the parallel computation power of today’s GPUs for time- and memory-efficiently solving large-scale IHPs under periodic boundary conditions with the density representation. 
However, it is challenging to make full use of the computing resources of GPU to realize the solver. The reasons are twofold.
First, since the solver contains multiple steps with different computational profiles, the choice of data structures and algorithms should be considered globally to be suitable for every step. 
Second, as the GPU memory is limited, the memory usage should be reduced to adapt to high resolution while ensuring accuracy and high efficiency.

\begin{figure*}[t]
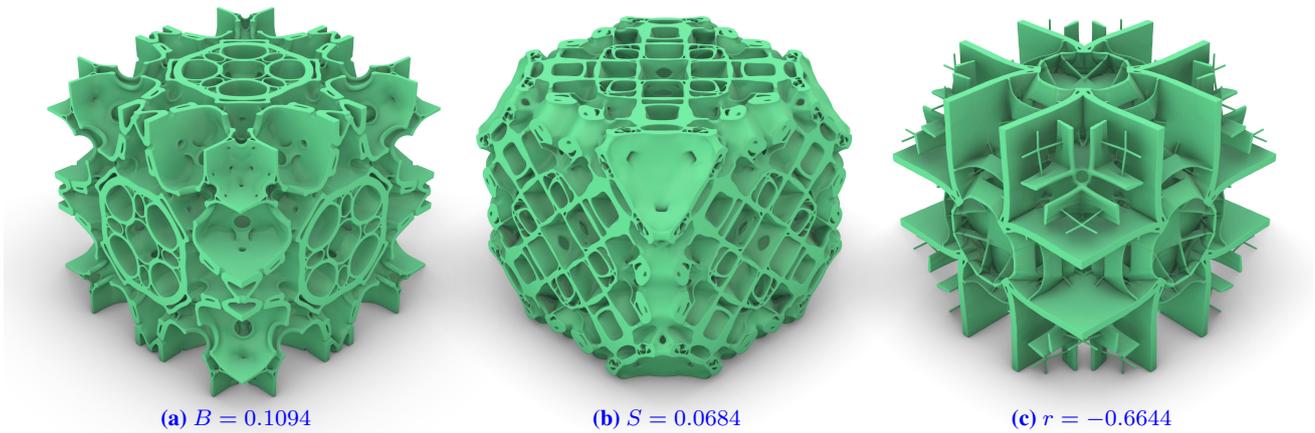

	\centering
	\begin{overpic}[width=0.99\linewidth]{teaser-rev-cc}
		{
			\put(12,1){\small \tb{\textbf{(a)} $B=0.1094$}}
			\put(45,1){\small \tb{\textbf{(b)} $S=0.0684$}}
			\put(77,1){\small \tb{\textbf{(c)} $r=-0.6644$}}
		}
	\end{overpic}
	\vspace{-2mm}
	\caption{
		Solving large-scale IHPs on the unit cell domain ($512\times 512\times 512$ elements) with three different objectives: bulk modulus $B$ (a), shear modulus $G$ (b), and Poisson's ratio $r$ (c). 
		\tb{The corresponding volume fractions are 0.3, 0.3 and 0.2, respectively.}
		We show the optimized objective values below models.
		%
	}
	\label{fig:teaser}
\end{figure*}

The goal of~\citep{wu2015system} is similar with ours, and they propose a high-performance multigrid solver with deep integration of GPU computing for solving compliance minimization problems.
However, since what they need to store is different from ours due to different problems, their data structure and multigrid solver are unsuitable for ours.
For example, we should store six displacement fields to evaluate the elastic matrix and perform sensitivity analysis, but this storage is a heavy burden for most GPUs. 
Besides, handling periodic boundary conditions is another difference.
%

This paper proposes an optimized, easy-to-use, open-source GPU solver for large-scale IHPs with periodic boundary conditions to design high-resolution 3D microstructures.
By exploring a software-level design space \tb{operating at only one GPU}, we present a favorable combination of data structures and algorithms to address the computational challenges of the desired solver.
Specifically, the mixed floating-point precision representation is deeply integrated into all components of the solver to achieve a favorable trade-off between memory usage, running time, and microstructure quality.
\tb{The mixed-precision formats under various precisions are tested, and we ultimately find that FP16/FP32 is the most suitable for IHPs. }
Besides, different types of memory are called properly and reasonably to significantly increase the number of optimizable finite elements on \tb{a} GPU.
We also provide test results to identify these favorable design choices.

We demonstrate the capability and superiority of our GPU solver by successfully optimizing the homogenized stiffness tensors, such as bulk modulus, shear modulus, and Poisson’s ratio, under the material volume-bounded constraint (Fig.~\ref{fig:teaser}).
%
\tb{In practice, our solver consumes less than \tb{40 seconds} for each iteration with a peak \tb{GPU} memory of 9 GB for high-resolution examples with $512^3 \approx 134.2$ million finite elements on an NVIDIA GeForce GTX \tb{1080Ti} GPU.}
Besides, we provide an easy-to-use framework for GPU implementation.
Specifically, the framework uses less than 20 lines of code to support various objective functions defined by the homogenization stiffness tensor.
Code for this paper is at \url{https://github.com/lavenklau/homo3d}.


\section{Related work}

\paragraph{Inverse homogenization problems.}
Solving IHPs~\citep{sigmund1994materials} to \tb{optimize the distribution of} materials is a powerful method to obtain superior mechanical properties under given load and boundary conditions.
Topology optimization is used to solve IHPs with different objectives, such as extreme shear or bulk moduli~\citep{gibiansky2000multiphase}, negative Poisson's ratios~\citep{theocaris1997negative,shan2015design,morvaridi2021hierarchical}, and extreme thermal expansion coefficients~\citep{sigmund1997design}.
Although there are several open-source codes for microstructure design~\citep{xia2015design,gao2021igatop}, most of them focus on 2D microstructure design.
We focus on developing an efficient GPU solver for large-scale IHPs to design 3D microstructures via the density-based method.

\paragraph{High-resolution topology problems.}
Several acceleration techniques for high-resolution topology problems are available, such as parallel computing~\citep{borrvall2001large,aage2015topology}, GPU computation~\citep{challis2014high}, adaptive mesh refinement~\citep{stainko2006adaptive,de2008topology,rong2022structural}.
Equipping and solving large-scale equilibrium equations is essential for slowing down the optimization process.
Thus, the geometric multigrid solver is used~\citep{briggs2000multigrid,zhu2010efficient,mcadams2011efficient,10.1111:cgf.14698}.
~\cite{wu2015system} present a high-performance multigrid solver with deep integration of GPU computing for compliance minimization problems.
PETSc~\citep{aage2015topology} is a large-scale topology optimization framework, where 
each iteration takes about 60s on 40 CPUs (240 cores) for a model with $288^3$ (23.8 million) elements. 
We focus on using the parallel computation power of today’s GPU and customize data structure and multigrid solver to solve large-scale IHPs of 3D microstructure design.
 \tb{A comparison with~\cite{aage2015topology} is shown in Fig.~\ref{fig:cmp-petsc}}

\paragraph{Mixed-precision methods.}
The IEEE standard provides for different levels of precision by varying the field width, e.g., \tb{16 bits (half precision)}, 32 bits (single precision) and 64 bits (double precision).
Double-precision arithmetic earns more accurate computations by suffering higher memory bandwidth and storage requirements.
However, \tb{half precision is 4 times speedup for a double precision~\cite{haidar2018harnessing}}
single-precision calculations take 2.5 times faster than the corresponding double-precision calculations~\citep{goddeke2010cyclic}.
Mixed-precision algorithms are proposed in many works for the trade-off between high efficiency and high precision~\citep{sun2008high,ben2020pop,zhang2019efficient,hosseini2023towards}.
\cite{liu2018narrow} developed a mixed-precision multigrid solver to accelerate the linear elasticity calculation.
We use the mixed-precision representation for a trade-off between memory usage, running time, and microstructure quality.

\section{Inverse Homogenization Problem}\label{sec:homo}

\subsection{Model}


The IHP is performed on a unit cell domain \tb{$\Omega$ = $[0,1]^3$}, which is evenly discretized into $M$ elements. %
Each element is assigned an density variable $\rho_e$ and a fixed volume $v_e$.
\tb{All density variables $\rho_e (e=1,\cdots,M)$ form a vector $\bm{\rho}$.}
IHP is formulated as follows:
\begin{equation} \label{eq:opt-prob}
	\begin{aligned}
		\underset{\bm{\rho}}{\min} \quad  &J = f(C^H(\bm{\rho})),
		\\
		\text{s.t.} \quad & \mathbf{K}\mathbf{u}=\mathbf{f},
		\\
		& \frac{\sum_{e=1}^M v_e \cdot \rho_e }{|\Omega|} \leq V,
		\\
		& \rho_{min}\leq \rho_{e} \leq 1, \quad \forall{e = 1,\cdots,M}.
	\end{aligned}
\end{equation}
Here, the objective $f(C^H(\bm{\rho}))$ is a function as the elastic matrix $C^H(\bm{\rho})$ to indicate mechanical properties.
The displacement field $\mathbf{u}$ is calculated by solving the equilibrium equation with six load cases $\mathbf{f}$ for three dimensions.
The stiffness matrix $\mathbf{K}$ is a function of the material properties in the elements.
$|\Omega|$ is the volume of the unit cell domain $\Omega$, $V$ is the prescribed volume fraction, and $\rho_{min} = 0.001$. 

\subsection{Homogenization}
Homogenization theory is typically used to determine the elastic tensor $E^H$ of 
a microstructure~\citep{shapehom2002}.
The derivation of the homogenized elasticity tensor involves a two-scale asymptotic expansion and boils down to solving the following \emph{cell problem} :
\begin{equation}
	\label{eq:cell-prob}
	\left\{
	\begin{aligned}
		\left.-\nabla \cdot\left(E:\left[\varepsilon\left(\mathbf{w}^{k l}\right)+e^{k l}\right]\right)\right)=0 \text { in } \Omega,\\
		\mathbf{w}^{kl}\left(\mathbf{x}\right)=\mathbf{w}^{kl}\left(\mathbf{x}+\mathbf{t}\right),\quad  \mathbf{x}\in\partial\Omega. \\
	\end{aligned}
	\right.
\end{equation}
Here, 
$e^{kl}$ with $k,l\in\{1,2,3\}$ is a unit tensor whose $(kl)$-th component equals to $1$ and other components equal to $0$.
The operator ``$:$'' means the double dot product of two tensors.
The second equation means $\mathbf{w}^{kl}$ is a periodic function whose period is $\mathbf{t}$ with $t_i=\pm 1, i=0,1,2$.
$\left(\nabla\boldsymbol{\chi}^{kl}\right)$ is the gradient of $\boldsymbol{\chi}^{kl}$  computed as $[\nabla\boldsymbol{\chi}^{kl}]_{ij}=\frac{\partial \chi^{kl}_i}{\partial x_j}$.
$E$ is the spatially varied elastic tensor of the base material. $\varepsilon(\mathbf{w}^{kl})=\frac{1}{2}\left(\nabla\mathbf{w}^{kl}+(\nabla\mathbf{w}^{kl})^\top\right)$ is the Cauchy strain tensor of the displacement field $\mathbf{w}^{kl}$.

After solving this problem, the homogenized elastic tensor is determined as:
\begin{equation}
	\label{eq:homtensor}
	E^H_{ijkl}=\frac{1}{|\Omega|}\int_\Omega (e^{ij}+\varepsilon(\mathbf{w}^{ij})):E:(e^{kl}+\varepsilon(\mathbf{w}^{kl}))\text{d}\Omega.
\end{equation}

For numerical computation, Finite Element Method is used to solve~\eqref{eq:cell-prob}~\citep{homnum2014}. 
We first enforce a macro strain on each element and compute the response force $\mathbf{f}_e^{ij}$:
\begin{equation}
	\label{eq:el-force}
	\mathbf{f}_e^{ij}=\mathbf{K}_e \boldsymbol{\chi}^{ij}_e,
\end{equation}
where $\boldsymbol{\chi}^{ij}_e$ is the displacement on element $e$'s vertices corresponding to the unit strain tensor $e^{ij}$.
The global force vector $\mathbf{f}^{ij}$ is assembled from element force vector $\mathbf{f}^{ij}_e$.
Based on the SIMP approach~\citep{bendsoe1989optimal}, $\mathbf{K}_e= \rho^p_e\mathbf{K}^0$ is the element stiffness matrix, where $\mathbf{K}^0$ is the element stiffness matrix of the element filled with base material and $p$ is a penalization factor.
The global stiffness matrix $\mathbf{K}$ is assembled from element stiffness matrix $\mathbf{K}_e$.
Then, we achieve the numerical solution $\mathbf{u}^{ij}$ by solving
\begin{equation} \label{eq:fem}
	\mathbf{K}\mathbf{u}^{ij}=\mathbf{f}^{ij}.
\end{equation}
After solving \eqref{eq:fem} for each pair of $ij$, the homogenized elastic tensor is computed as
\begin{equation}\label{eq:elasttensor}
	E^H_{ijkl}=\frac{1}{|\Omega|} \sum_e (\boldsymbol{\chi}^{ij}_e-\mathbf{u}^{ij}_e)^\top\mathbf{K}_e(\boldsymbol{\chi}^{kl}_e-\mathbf{u}^{kl}_e),
\end{equation}
where $i,j,k,l\in \{1,2,3\}$ and $\mathbf{u}_e^{ij}$ means the components of $\mathbf{u}^{ij}$ on the element $e$.
Using the engineering notation with $11\rightarrow 0$, $22\rightarrow 1$, $33\rightarrow 2$, $12\rightarrow 3$ , $23\rightarrow 4$  and $13\rightarrow 5$ , the elasticity tensor, i.e., $E_{ijkl}^H$ in~\eqref{eq:elasttensor}, is rewritten as
\begin{equation}\label{eq:elasttensor}
	C^H_{ij}=\frac{1}{|\Omega|} \sum_e (\boldsymbol{\chi}^{i}_e-\mathbf{u}^{i}_e)^\top\mathbf{K}_e(\boldsymbol{\chi}^{j}_e-\mathbf{u}^{j}_e).
\end{equation}
The objective $f(C^H)$ is a user-defined function as the components of $C^H$.
Its gradient is computed as:
\begin{equation}\label{eq:obj-grad}
	\frac{\partial f}{\partial \rho_e}=\sum_{ij}\frac{\partial f}{\partial C^H_{ij}}\frac{\partial C^H_{ij}}{\partial \rho_e},
\end{equation}
where 
\begin{equation}\label{eq:ch-sens}
	\frac{\partial C^H_{ij}}{\partial \rho_e}=\frac{1}{|\Omega|}\sum_e p \rho_e^{p-1}
	(\boldsymbol{\chi}_e^{i}-\mathbf{u}_e^{i})^\top\mathbf{K}^0
	(\boldsymbol{\chi}_e^{j}-\mathbf{u}_e^{j}).
\end{equation}

\subsection{Optimization model}
\paragraph{Solver for IHP~\eqref{eq:opt-prob}.}
We solve IHP in an iterative manner.
In each iteration, the following four steps are performed:
\begin{enumerate}
	\item Compute the displacement field $\mathbf{u}$ by solving~\eqref{eq:fem}. 
	\item Compute the homogenized elastic tensor $C^H$ via~\eqref{eq:elasttensor} and the objective function $f(C^H)$.
	\item Perform sensitivity analysis, i.e., evaluate the gradient $\frac{\partial f}{\partial \bm{\rho}}$ via~\eqref{eq:obj-grad} and~\eqref{eq:ch-sens}
	\item Update density $\boldsymbol{\rho}$ using $\frac{\partial f}{\partial \bm{\rho}}$ based on the Optimal Criteria (OC) method~\citep{sigmund200199}.  
\end{enumerate}


\paragraph{Multigrid solver.}
Solving~\eqref{eq:fem} to compute $\mathbf{u}$ for large-scale problems is time-consuming and memory-intensive.
To reduce the time and memory overhead, the multigrid solver is used~\citep{DICK2011801,wu2015system,liu2018narrow}.
The main idea of multigrid is to solve a coarse problem by a global correction of the fine grid solution approximation from time to time to accelerate the convergence.

Our first level grid is the cell domain $\Omega$.
Then, we recursively divide $\Omega$ to construct a hierarchy of coarse grids.
To transfer data between grids of levels $l$ and $l+1$, we use trilinear interpolation and its transpose as the restriction operator, denoted as $I^{l}_{l+1}$ and $R^{l+1}_{l}$, respectively.
Based on Galerkin rule, the numerical stencil on the level $l+1$ is determined from the level $l$ as
$\mathbf{K}^{l+1}=R^{l+1}_{l}\mathbf{K}^{l}I^{l}_{l+1}$.
Then, the V-cycle of the multigrid solver, with the Gauss-Seidel relaxation as the smoother, is employed to effectively decrease the residual on the first level grid until convergence.
In addition, 
\tb{the coarsening process would be stopped if the subsequent coarsened mesh becomes smaller than $4 \times 4 \times 4$ while using coarsening ratio of 1 : 2.}
The pseudocode of the V-cycle is outlined in Alg.~\ref{alg:multigrid}.


\IncMargin{0.5em}
\begin{algorithm}[t]
	\caption{V-cycle in multigrid solver}\label{alg:multigrid}
	\SetCommentSty{mycommfont}
	\SetKwInOut{AlgoInput}{Input}
	\SetKwInOut{AlgoOutput}{Output}
	\For{$l=0,\cdots,L-1$}
	{
		\If{$l>0$}{$\mathbf{u}^{l}\gets 0$}
		Relax $\mathbf{K}^{l}\mathbf{u}^{l}=\mathbf{f}^{l}$;\tcp*[f]{Relaxation}\\
		$\mathbf{r}^{l}=\mathbf{f}^{l}-\mathbf{K}^{l}\mathbf{u}^{l}$;\tcp*[f]{Residual update}\\
		$\mathbf{f}^{l+1}=R^{l+1}_{l}\mathbf{r}^{l}$;\tcp*[f]{Restrict residual}\\
	}
	Solve $\mathbf{K}^{L}\mathbf{u}^{L}=\mathbf{f}^{L}$ directly;\tcp*[f]{Solve on coarsest level}\\
	\For(\tcp*[f]{Go up in the V-cycle}){$l=L-1,\cdots,0$}
	{
		$\mathbf{u}^{l}\gets \mathbf{u}^{l}+I_{l+1}^{l}\mathbf{u}^{l+1}$;\tcp*[f]{Interpolate error \& correct}\\
		Relax $\mathbf{K}^{l}\mathbf{u}^{l}=\mathbf{f}^{l};$\tcp*[f]{Relaxation}\\
	}
\end{algorithm}
\DecMargin{0.5em}

\section{Optimized GPU Scheme for solving IHPs}\label{sec:method}




We describe our optimized GPU scheme for solving large-scale IHPs using a GPU-tailored data structure (Section~\ref{sec:data}), a dedicated multigrid solver (Section~\ref{sec:mg}), and an efficient evaluation of the elastic matrix and sensitivity (Section~\ref{sec:elmat-sens}).


\subsection{Data structure tailored to solve IHPs}\label{sec:data}
\paragraph{Data for each vertex.}
For each vertex $\mathbf{v}$ of each level's mesh, we store the numerical stencil $\mathbf{K}_v$, the displacement $\mathbf{u}_v$, the force $\mathbf{f}_v$, and the residual $\mathbf{r}_v$ in the multigrid implementation.
The numerical stencil $\mathbf{K}_v$ consists of 27 matrices of dimension $3\times3$, each of which corresponds to one adjacent vertex.
The displacement, force, and residual are all vectors with three components, named \emph{nodal vectors}.
From the perspective of implementation, a numerical stencil can be regarded as composed of $27 \times 3$ nodal vectors. 
On the first level mesh, the density variable is stored for each element, which is used to assemble the numerical stencil on the fly.
To handle different boundary conditions and facilitate the Gauss-Seidel relaxation, a 2-byte flag is stored for each vertex and element on each level's mesh, named \emph{vertex flag} and \emph{element flag}, respectively.

\paragraph{Mixed floating-point precision representations.}

\tb{By employing mixed floating-point precision, a balance between computational accuracy and performance can be achieved.
Since a large-scale IHP involves large elements calculation and requires significant memory usage, applying mixed precision can help conserve memory and enable more efficient data storage.
The throughput of lower precision format is usually much higher than that of higher precision. Only using the higher precision through the computation leaves the lower precision pipeline underutilized, wasting significant computing resources.
Through extensive testing (Table~\ref{tab:pre}), we are pleasantly surprised to find that using mixed floating-point precision (FP16/FP32) can maximize computational efficiency and memory storage under the tolerance relative residual $r_{rel} = 10^{-2}$, computed as $\|\mathbf{r}\|_2/\|\mathbf{f}\|_2$ (see more discussions about $r_{rel}$ in Fig.~\ref{fig:relaxation-residualupdate} and Section~\ref{sec:mg}). 
The numerical stencils are stored in half-precision (FP16), and the rest vectors are stored in single-precision (FP32).}



\paragraph{Memory layouts.}
\tb{Nodal vectors are all stored in the Structure of Array (SoA) format.}
Namely, the same component of stencil or nodal vector of all vertices is stored together in an array, and different components are stored in different arrays.
\tb{The numerical stencils are stored in Array of Structure (AoS) format.}

Specifically, the eight-color Gauss-Seidel relaxation is used for the parallelization, which partitions the vertices into eight independent subsets and parallelizes the computation within each subset.
Our memory layout should provide an efficient memory access pattern for this procedure.
Given a vertex with an integer coordinate $(x_0,x_1,x_2)$, it belongs to the subset $(x_0\,\text{mod}\, 2)2^0 + (x_1\,\text{mod}\, 2) 2^1 + (x_2\,\text{mod}\,2) 2^2$.
To exploit the high memory bandwidth and leverage the coalesced memory transaction, we use a similar memory layout as~\citep{DICK2011801}, i.e., the data on the vertices of the same subset are grouped and different subsets are stored in the consecutive memory block.


In most modern GPUs, only the  32-, 64-, or 128-byte segments of device memory that are aligned to multiple of their size can be read or written by memory transactions. To maximize the memory throughput, data should be organized in such a way that the i-th thread of a warp (containing 32 threads) accesses the i-th 32-bit (or 64-bit) word of a 128-byte segment at single (or double) precision.
Thus, we add a few ``ghost'' vertices to supplement each subset so that its number of vertices is a multiple of 32.
Then, we can assign one warp to each group of 32 vertices for each subset.
Each warp accesses an aligned 128-byte (or 256-byte) segment for single-precision (or double-precision) nodal vectors. 
Ideally, each warp's memory access to the same adjacent vertices of the group of 32 vertices is also coalesced as these adjacent vertices are consecutive in another subset's memory block.

\begin{figure}[t]
	\centering
	\begin{overpic}[width=0.9\linewidth]{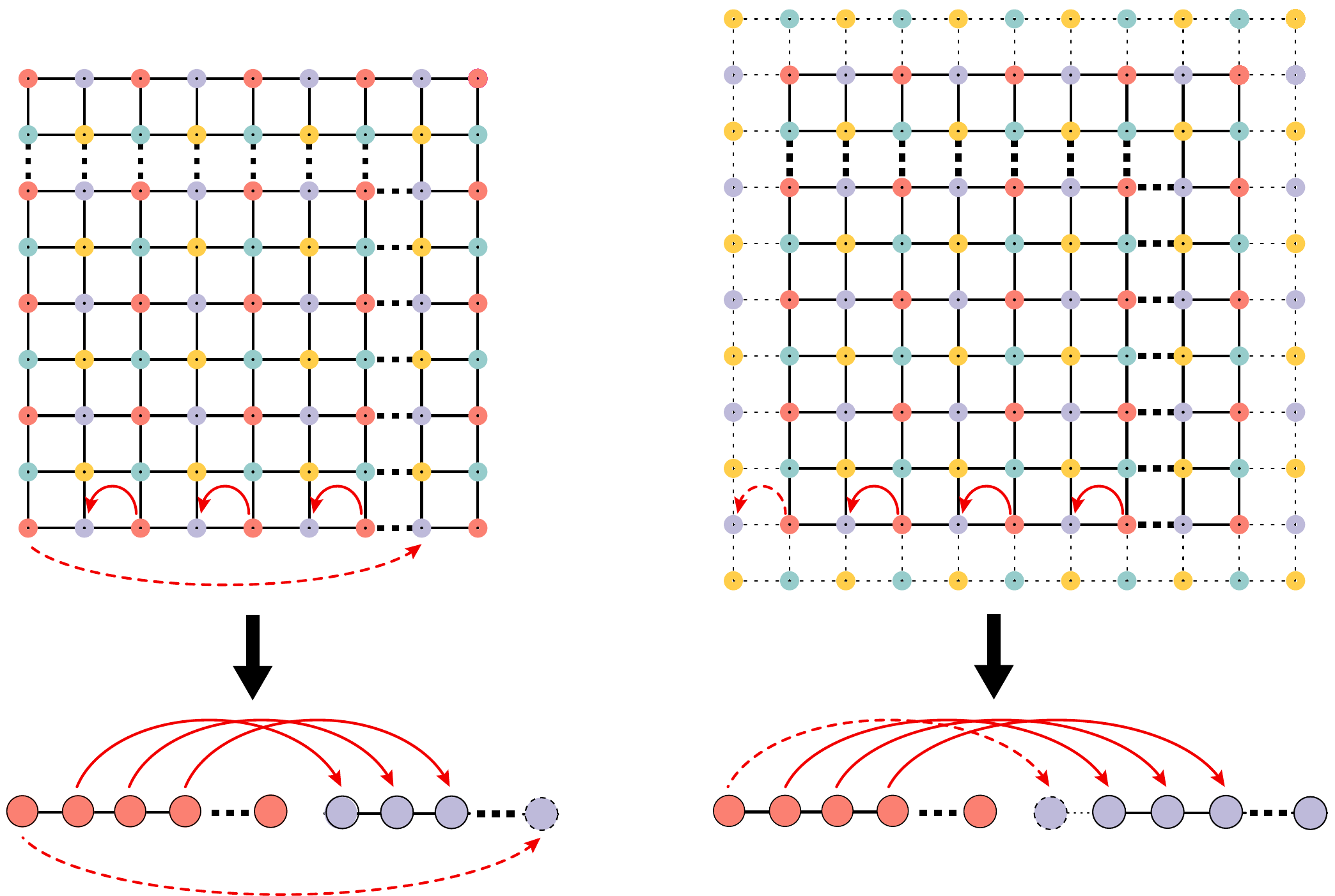}
		{
		}
	\end{overpic}
	\vspace{-1mm}
	\caption{
		2D illustration for periodic boundary conditions.
		In the Gauss-Seidel relaxation, the vertices are partitioned into four subsets.
		Each subset is shown in one color.
		We use the memory access to the left adjacent vertex (the red solid or dotted arrow) for illustration.
		For the red dotted arrow, memory access starts from the vertex of the left boundary.
		Before padding (upper left), memory access from the left boundary vertex is discontinuous (lower left), leading to uncoloasced memory transactions.
		After padding (upper right), such memory access becomes continuous (lower right); thus, the memory access is coloasced.
	}
	\label{fig:padmesh}
\end{figure}

\paragraph{Padding layers for periodic boundary conditions.}
The coarse mesh in the multigrid solver should inherit the periodic feature of the fine mesh.
Due to periodic boundary conditions, the vertices on the opposite boundaries are the same.
Hence, when restricting residuals or numerical stencils from the fine grid to the coarse grid, the vertices on one side of the boundary should add the transferred data from the neighbor on the other side.

In our GPU implementation, we pad a layer of vertices and elements around the mesh (Fig.~\ref{fig:padmesh}).
Those padded vertices and elements are copies of their periodic equivalents and are updated when their copied vertices or elements change.
After the padding, the memory layout is updated to incorporate the padded vertices.
Then, the restriction can transfer data from neighbors regardless of periodic boundary conditions.
We do not execute computations on the padded vertices or elements and they only provide data to their neighbors.
This padding leads to a more efficient memory access pattern (see an example in Fig.~\ref{fig:padmesh}).

\paragraph{Accessing data in memory.}
Our grid is highly regular as it is evenly divided from a cube.
Given the integer coordinate $(x_0,x_1,x_2)$ of a vertex, the memory location is:
\begin{equation}\label{eq:p-to-mem}
	p = p_\text{base}^\text{Id} + \lfloor x_0 / 2 \rfloor  + \left(\lfloor x_1 / 2 \rfloor  + \lfloor x_2 / 2 \rfloor N_1^\text{Id}\right) N_0^\text{Id},
\end{equation}
where $\text{Id} = (x_0\,\text{mod}\, 2)2^0 + (x_1\,\text{mod}\, 2) 2^1 + (x_2\,\text{mod}\,2) 2^2$ is the index of the subset, $p_\text{base}^\text{Id}$ denotes the start address of the memory block of the subset $\text{Id}$, $N^\text{Id}_i,i=0,1,2$ is the number of vertices of the subset $\text{Id}$ along three axes:
$$N^\text{Id}_i=\lfloor \left(N_i - O^\text{Id}_i\right)/2\rfloor + 1,$$
where $N_i,i=0,1,2$ is the number of elements along three axes, $O^\text{Id}_i$ is the origin of the subset \text{Id} defined as $O^\text{Id}_0=\text{Id} \,\text{mod}\,  2,O^\text{Id}_1=\lfloor \text{Id} / 2\rfloor \,\text{mod}\, 2,O^\text{Id}_2 =\lfloor \text{Id} / 4\rfloor$.

Since the coordinates of one vertex's adjacent elements or vertices can be calculated by offsetting the position of itself, we can compute their memory locations easily.
Hence, we do not store the topology information, e.g., the index of the adjacent vertex, to reduce a large amount of memory.

\paragraph{Different GPU memory types.}
In high-resolution problems, the storage for the nodal vectors, e.g., nodal displacements, is enormous.
For example, \tb{we need about 1.5GB memory for one single-precision nodal displacement field of a grid with a resolution of $512^3$.}
Six displacement fields should be stored to evaluate the elastic matrix and perform sensitivity analysis.
Besides, they are proper initializations for solving~\eqref{eq:fem} in the next iteration.
However, it costs nearly \tb{9GB} of memory, which is unaffordable for most GPUs.
Thus, our GPU implementation stores them in unified memory that supports oversubscription.
Hosts and devices can access unified memory and the CUDA underlying system manages its physical location.
We provide methods to reduce the performance loss of unified memory in Section~\ref{sec:elmat-sens}.
%

Except for these displacement fields, the storage for the nodal vectors of the multigrid solver is allocated and resident on the device memory, which is the same as the density variable, vertex, and element flag.
The numerical stencils are also stored in the device memory except for the first level mesh, where we assemble the numerical stencil using densities on the fly.
The frequently used data are cached on constant memory, such as the template matrix, pointers to vertices, element data, and grid information like resolution.

\begin{figure}[t]
	\centering
	\begin{overpic}[width=0.99\linewidth]{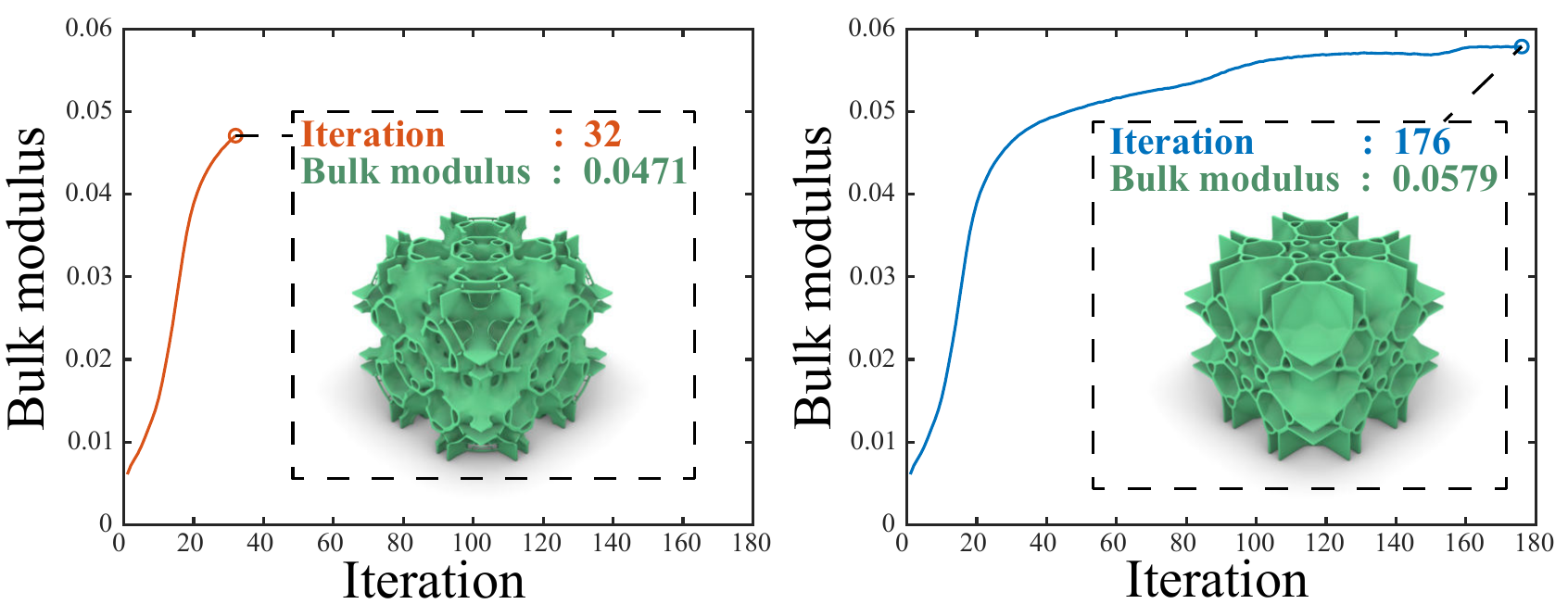}
		{
			\put(20,-3.5) {\small \textbf{(a)} W/O}
			\put(70,-3.5) {\small \textbf{(b)} With}
		}
	\end{overpic}
	\vspace{2mm}
	\caption{
		\tb{
		An ablation study of the Dirichlet boundary for maximizing the bulk modulus with the resolution of $128^3$ and the volume ratio of $0.2$.
	}
	}
	\label{fig:mix-fp32-err}
\end{figure}

\begin{figure}[t]
	\centering
	\begin{overpic}[width=0.99\linewidth]{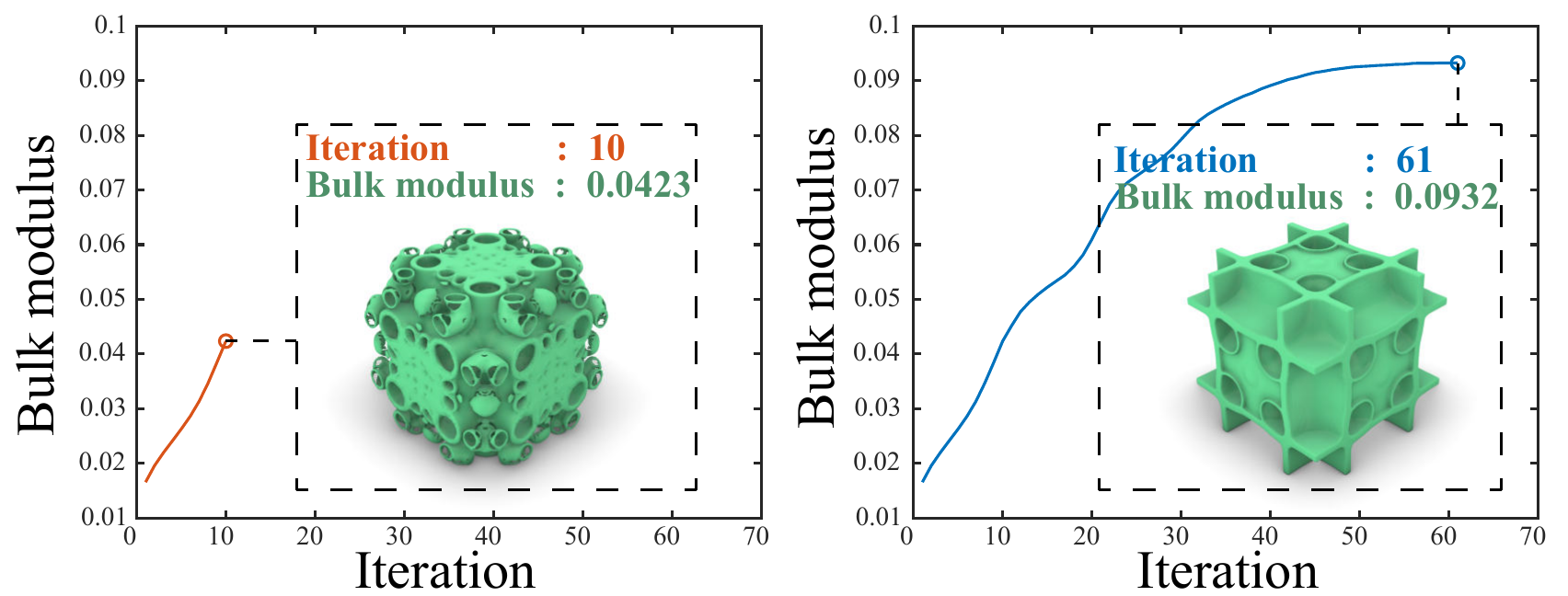}
		{
			\put(13,-3.5) {\small \textbf{(a)} W/O removing}
			\put(63,-3.5) {\small \textbf{(b)} With removing}
		}
	\end{overpic}
	\vspace{2mm}
	\caption{
		\tb{
			Singular stiffness matrices.
			We optimize the bulk modulus, and the initial density fields of (a) and (b) are the same.
			If our multigrid solver does not remove the component belonging to the null space of numerical stencil, it diverges at the 10-st iteration (a).
		}
	}
	\label{fig:rigid}
\end{figure}

\subsection{Dedicated multigrid solver}\label{sec:mg}

\paragraph{Singular stiffness matrices.}
Due to the loss of precision caused by the mixed-precision scheme and the high resolutions, the multigrid solver may diverge with a numerical explosion.
\tb{We find in practice that these situations may be caused by (1) insufficient Dirichlet boundary conditions and (2) no materials at corners during optimization.}

During homogenization, the eight corner vertices of the unit cube domain are usually selected as the fixed vertices.
\tb{This amounts to adding Dirichlet boundary conditions at the eight corner vertices to the cell problem~\eqref{eq:cell-prob} to guarantee a unique solution (see Fig.~\ref{fig:mix-fp32-err}~(b)).
Otherwise, the global stiffness matrix becomes singular (see Fig.~\ref{fig:mix-fp32-err}~(a))}.

However, as the density field evolves during the optimization, it often tends to be zero near the corners. 
Accordingly, the solid part gets isolated from the corners.
Then, the global stiffness matrix is again becoming singular.
To handle such a problem, we remove the component belonging to the numerical stencil's null space from the restricted residual before solving the system on the coarsest mesh, similar to~\citep{elasticTexture,10.1111:cgf.14698}.
\tb{We show an example with the resolution $128^3$ under the volume fraction 0.3 in Fig.~\ref{fig:rigid}.}

\paragraph{Enforce macro strain.}
The response force on the vertex $\mathbf{v}$ from an enforced macro strain $e^{kl}$ is:
\begin{equation}
	\mathbf{f}^{i}_{v} = \sum_{e=0}^7\rho_e^p\left(\sum_{v_j=0}^7\mathbf{K}^0_{[7-e,v_j]}\boldsymbol{\chi}^{i}_{e,v_j}\right),
\end{equation}
where $\mathbf{K}^0_{[i,j]}$ denotes $(i, j)$-th $3\times3$-block of $\mathbf{K}_0$, and $\boldsymbol{\chi}^{i}_{e,v_j}$ is the displacement on vertex $\mathbf{v}_j$ for the macro strain, where the superscript $i$ is the engineering notation for $kl$.

To enforce macro strain $e^{kl}$, we assign one thread for each vertex.
Each thread traverses the incident elements of its assigned vertex and accumulates the response force of each element on this vertex.
On the vertex $\mathbf{v}_j$, we have
\begin{equation}
	\label{eq:macro-strain}
	\begin{aligned}
		\boldsymbol{\chi}^{0}_{e,v_j}&=(x^{v_j}_0,0,0)^\top,\quad  &\boldsymbol{\chi}^{3}_{e,v_j}&=(x^{v_j}_1/2,x^{v_j}_0/2,0)^\top,\\
		\boldsymbol{\chi}^{1}_{e,v_j}&=(0,x^{v_j}_1,0)^\top,\quad  &\boldsymbol{\chi}^{4}_{e,v_j}&=(0,x^{v_j}_2/2,x^{v_j}_1/2)^\top,\\
		\boldsymbol{\chi}^{2}_{e,v_j}&=(0,0,x^{v_j}_2)^\top,\quad  &\boldsymbol{\chi}^{5}_{e,v_j}&=(x^{v_j}_2/2,0,x^{v_j}_0/2)^\top,
	\end{aligned}
\end{equation}
where $\mathbf{x}^{v_j} = (x^{v_j}_0,x^{v_j}_1,x^{v_j}_2)$ is the coordinate of the vertex $\mathbf{v}_j$.
Due to the accuracy loss of the half-precision stiffness matrix, the translation of the nodal displacement causes a response force that cannot be ignored numerically.
Consequently, the absolute position of $\mathbf{v}_j$ affects the response force.
Hence, we use the relative coordinate of $\mathbf{v}_j$ in the element rather than the coordinate in the entire grid as $\mathbf{x}^{v_j}$.

\paragraph{Relaxation and residual update.}
To implement the eight-color Gauss-Seidel relaxation, we serially launch one computation kernel for each subset of the vertices.
The performance bottleneck of the multigrid solver is the Gauss-Seidel relaxation and residual update on the first level mesh.
Central to both procedures is to compute $\mathbf{Ku}$ on each vertex $\mathbf{v}$:
\begin{equation}\label{eq:st-assemb}
	[\mathbf{Ku}]_v = \sum_{e=0}^7\rho_e^p\left(\sum_{v_j=0}^7\mathbf{K}^0_{[7-e,v_j]}\mathbf{u}^{v_j,e}\right),
\end{equation}
where 
the subscript $e$ is the incident element of the vertex $\mathbf{v}$, and $\mathbf{u}^{v_j,e}$ is the nodal displacement on $v_j$-th vertex of element $e$.
The residual is then updated as
\begin{equation}
	\mathbf{r}_v = \mathbf{f}_v - [\mathbf{Ku}]_v.
\end{equation}
We introduce two notations for the relaxation:
\begin{equation}
	\label{eq:st-assemb-gs}
	\begin{aligned}
		&\mathbf{M}_v = \sum_{e=0}^7\rho_e^p\left(\sum_{v_j=0,v_j\ne7-e}^7\mathbf{K}^0_{[7-e,v_j]}\mathbf{u}^{v_j,e}\right), \\
		&\mathbf{S}_v = \sum_{e=0}^7\rho_e^p \mathbf{K}^0_{[7-e,7-e]},
	\end{aligned}
\end{equation}
where we use $\mathbf{M}_v$ to denote the modified $[\mathbf{Ku}]_v$ and $\mathbf{S}_v$ to denote the sum of $3\times3$ diagonal block of $\mathbf{K}^0$.
Then, the Gauss-Seidel relaxation is performed via the following linear system to update $\mathbf{u}_v$:
\begin{equation}\label{eq:st-gs-relx}
	\mathbf{S}_v\mathbf{u}=\mathbf{f}_v - \mathbf{M}_v.
\end{equation}

To perform these computations via GPU, we first dispatch eight warps for each group of 32 consecutive vertices.
Each warp accumulates the contribution of one incident element in~\eqref{eq:st-assemb} or \eqref{eq:st-assemb-gs}.
Then, we compute the total sum by a block reduction.

\begin{figure}[t]
	\centering
	\begin{overpic}[width=0.75\linewidth]{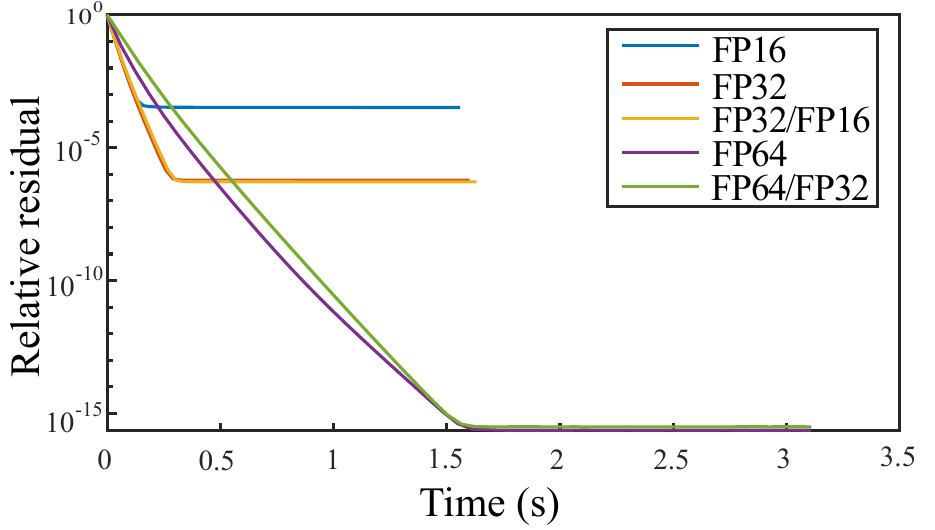}
			{
				}
		\end{overpic}
	\vspace{-2mm}
	\caption{
			\tr{
		}
		\tb{
The temporal evolution of the relative residual (i.e., $\|\mathbf{r}\|_2/\|\mathbf{f}\|_2$) in a multigrid solver with different precisions (FP16, FP32, FP64, FP32/FP16 and FP64/FP32).
We generate a random density between 0 and 1 for each element to form a random density field.
 The solver was iterated for 50 V-cycles.
	}
		}
	\label{fig:relaxation-residualupdate}
\vspace{2mm}
\end{figure}

\begin{figure}[t]
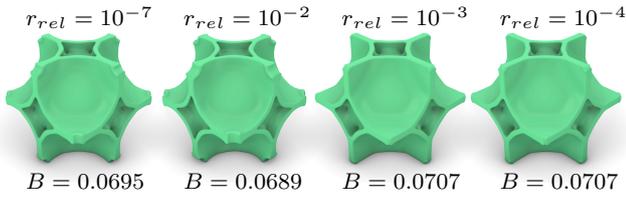

	\centering
	\begin{overpic}[width=0.99\linewidth]{diffrel_cc}
		{
			\put(5,25){\small $r_{rel}=10^{-7}$}
			\put(30,25){\small $r_{rel}=10^{-2}$}
			\put(55,25){\small $r_{rel}=10^{-3}$}
			\put(80,25){\small $r_{rel}=10^{-4}$}
			\put(5,-1){\small $B=0.0695$}
			\put(30,-1){\small $B=0.0689$}
			\put(55,-1){\small $B=0.0707$}
			\put(80,-1){\small $B=0.0707$}
		}
	\end{overpic}
	\vspace{0mm}
	\caption{
	\tb{
	The optimization results of the bulk modulus $B$ with different relative residuals $r_{rel}$s. The leftmost result is the baseline achieved with FP64, while the remaining three results are obtained using a mixed-precision scheme combining FP32 and FP16. 
	}
	}
	\label{fig:femrelthres}
\end{figure}

\tb{
To use the computational power of modern GPUs, lower precision representations such as FP32 and FP16 are preferred over FP64 due to their higher throughput and smaller bandwidth requirement. However, a trade-off exists between computational efficiency and accuracy.}
\tb{We explore different combination schemes of precision representations within our multigrid solver (Fig.~\ref{fig:relaxation-residualupdate} and Table~\ref{tab:pre}).
It is observed that  compared to the single-precision scheme, the mixed-precision scheme achieves a comparable relative residual using less memory.
Based on these comparisons, we have identified that the combination of FP32 and FP16 yields the best results within the specified tolerance error.}

\tb{We also test different $r_{rel}$s of the equilibrium equation in Fig.~\ref{fig:femrelthres}.
Again, the structures are almost the same, and the differences in bulk modulus are less than  3\%. }

\paragraph{Restriction and prolongation.}
We follow~\citep{DICK2011801} to restrict residuals and prolong displacements, except that the index of the adjacent vertex is computed via~\eqref{eq:p-to-mem} instead of being loaded from global memory.


\paragraph{Assembling numerical stencils for coarse grids.}
Since the numerical stencil on the first level is not stored, we assemble the numerical stencil on the second level as follows:
\begin{equation}\label{eq:st-level-2}
	[\mathbf{K}_v]_{v_k} = \sum_{e\in N(v)}\sum_{v_i=0}^{7}w^{v}_{e,v_i}\sum_{v_j=0}^7w^{v_k}_{e,v_j} \rho_e^p\mathbf{K}^0_{[v_i,v_j]}.
\end{equation}
Here $[\mathbf{K}_v]_{v_k}\in\mathbb{R}^{3\times3}$ is the numerical stencil of the vertex $\mathbf{v}$ to its adjacent vertex $\mathbf{v}_k$ on the second level grid.
$N(v)$ is the set of elements on the first level grid covered by the adjacent elements of $\mathbf{v}$ on the second level.
The weight $w^v_{e,v_i}$ is:
\begin{equation}\label{eq:st-assemb-w}
	w^v_{e,v_i}=\prod_{k=1}^3 \frac{d-|x^{v}_k-x^{e,v_i}_k|}{d},
\end{equation}
where $(x^{v}_0,x^{v}_1,x^{v}_2)$ and $(x^{e,v_i}_0,x^{e,v_i}_1,x^{e,v_i}_2)$ are the coordinates of $\mathbf{v}$ and the $v_i$-th vertex of the element $e$, respectively, and $d$ is the length of element on the second level. The weight $w^{v_k}_{e,v_j}$ is defined in the same way.

\tb{
In the GPU implementation, we assign one thread to each vertex in the second level.
Each thread iterates through its $27$ neighboring vertices and accumulates the summands in~\eqref{eq:st-level-2}.
More specifically, in each loop, the thread accesses the density value in~\eqref{eq:st-level-2} from global memory, computes its power, and then multiplies it by the weights and the $3\times3$ block $\mathbf{K}^0_{[v_i,v_j]}$. The resulting product is then summed into a $3\times3$ matrix in local memory, which is written back to global memory at the end of each loop.
It is worth noting that, due to the presence of many zero weights in~\eqref{eq:st-level-2}, we only perform computations on the first-level elements that are covered by both the incident elements of the assigned vertex and the current looping neighboring vertices in the second level.
 This excludes many summands with zero weight.
}
If the resolution of the second level grid is high, e.g., the number of elements exceeds $128^3$, the memory cost is reduced by the non-dyadic coarsening strategy~\citep{wu2015system} that directly transfers the numerical stencil from the first level to the third level.

The numerical stencil on the higher-level grid (e.g., third level, fourth level, etc) is assembled in the same way.
Specifically, we use the numerical stencil on the fine grid to assemble the numerical stencil on the coarse grid, e.g., second level for third level.
First, our GPU implementation assigns a thread for each vertex on the coarse grid.
The thread loops 9 times to compute the 9 entries for all $3\times3$ matrices of the numerical stencil.
In each loop, the thread loads one of the 9 entries from the numerical stencil of its adjacent vertices on the fine grid and computes the weighted sum.
Then, the thread writes the sum back to the global memory.


\subsection{Elastic matrix evaluation and sensitivity analysis.}\label{sec:elmat-sens}
\begin{figure}[t]
	\centering
	\begin{overpic}[width=0.9\linewidth]{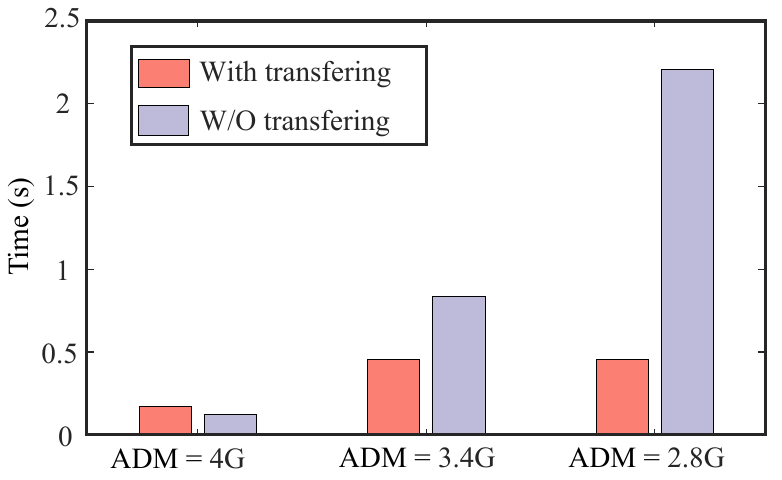}
		{
		}
	\end{overpic}
	\vspace{-2mm}
	\caption{
	\tb{	An ablation study of transferring the displacement fields to device memory.
		We report the timings for performing sensitivity analysis on a $256^3$ grid with different available device memories (ADMs).} 
		%
	}
	\label{fig:elastic-sensitivity}
\end{figure}

\paragraph{Handling unified memory.}
Evaluating elastic matrix and sensitivity heavily depends on the six displacement fields stored in the unified memory.
The performance loss of the unified memory increases the time cost for both operations.
We find in practice that the \tb{FP32 precision} displacement is necessary for numerical stability when solving FEM, whereas the \tb{FP16 precision} is enough to evaluate the elastic matrix and perform sensitivity analysis.
To reduce such performance loss, we first launch a kernel to cast the \tb{FP32 precision} displacements to \tb{FP16 precision} and then store them in the memory of the displacement, the residual, and the force on the first level mesh (Fig.~\ref{fig:elastic-sensitivity}).

\paragraph{Evaluation.}
The elastic matrix and sensitivity are evaluated similarly.
Eight warps are assigned for each group of 32 consecutive elements.
The first 6 warps compute $\boldsymbol{\chi}_e^{i}-\mathbf{u}_e^{i}, i\in\{0,1,2,3,4,5\}$ and store it in the shared memory, where $\boldsymbol{\chi}_e^{i}$ is computed on the fly and $\mathbf{u}_e^{i}$ is loaded from the memory.

To evaluate the elastic matrix in our GPU implementation, we first split $\left(\boldsymbol{\chi}_e^{i}-\mathbf{u}_e^{i}\right)^\top\mathbf{K}_e\left(\boldsymbol{\chi}_e^{j}-\mathbf{u}_e^{j}\right)$ in~\eqref{eq:elasttensor} into 8 summands by dividing $\mathbf{K}_e$ into eight 3-row blocks, which are dispatched into the eight warps, respectively.
Then, a block reduction is performed to get the product in the first warp for each element.
Since we aim to sum over all elements according to~\eqref{eq:elasttensor}, we do a warp reduction using the warp shuffle operation to compute the sum over the 32 elements in the first warp before writing it to memory.
Finally, several parallel reduction kernels are serially launched to compute the sum of the results produced by the last step.

For sensitivity analysis, the split product becomes
$\rho_e^{p-1}\frac{\partial f}{\partial C^{H}_{ij}}\left(\boldsymbol{\chi}_e^{i}-\mathbf{u}_e^{i}\right)^\top\mathbf{K}^0\left(\boldsymbol{\chi}_e^{j}-\mathbf{u}_e^{j}\right)$
with a constant coefficient.
We do not sum over all the elements as the sensitivity is computed for each element according to~\eqref{eq:ch-sens}.
\section{An easy-to-use framework}
\label{sec:autodiff}

\subsection{Setup}
Users can clone this framework or fork the current master branch from the GitHub repository (\url{https://github.com/lavenklau/homo3d}).
The compilation and runtime environment mainly requires \tb{CUDA 11, gflags, Eigen3, glm, and OpenVDB}.
The main classes are listed and explained in the supplementary material.

\subsection{Compiling and code invoking}\label{sec:compile-invoke}
The framework provides a good user interface.
After installing the framework, the user can use the following steps to design 3D microstructure:
\begin{itemize}
	\item The initialization includes optimization parameters, the design domain, and its resolution:
	\begin{lstlisting}
Homogenization hom(config);
TensorVar<float> rho(config.reso[0],config.reso[1],config.reso[2]);
initDensity(rho, config);
	\end{lstlisting}
	where \lstinline{config} is a parsed configuration file with command line arguments, including the Young's modulus and Poisson's ratio of the base material, grid resolution, volume fraction, initialization type, symmetry requirement, etc.
	\item Define the material interpolation method based on the SIMP approach~\citep{bendsoe1989optimal}:
	\begin{lstlisting}
auto rhop = rho.conv(radial_convker_t<float,Spline4>(config.filterRadius)).pow(3);
	\end{lstlisting}
	where \lstinline{conv(radial_convker_t<float, Spline4>(1.2))} means a convolution operation with the kernel \lstinline{radial_convker_t<float, Spline4>(1.2)}, which is same with the filtering method of~\citep{wu2015system}.
	\tb{The periodic filter kernel is discussed in Section~\ref{sec:diff-objs}.}
	\item Create an elastic matrix from the design domain \lstinline{hom} and the material  interpolation method \lstinline{rhop}:%
	\begin{lstlisting}
auto Ch = genCH(hom, rhop);
	\end{lstlisting}
	\item Define the objective function $f(C^H)$, e.g., the following objective is to maximize the bulk modulus:
	\begin{equation}\small
		\label{eq:bulk-obj}
		f(C^H)=-\frac{1}{9}\left(
		C^H_{00}+C^H_{11}+C^H_{22}+2(C^H_{01}+C^H_{02}+C^H_{12})\right).
	\end{equation}
	The code is written as
	\begin{lstlisting}
auto objective = -(Ch(0, 0) + Ch(1, 1) + Ch(2, 2) + (Ch(0, 1) + Ch(0, 2) + Ch(1, 2)) * 2) / 9.f;
	\end{lstlisting}
	\item Define the optimization process. We create an optimizer and begin the main optimization loop.
	In each iteration, we evaluate the objective, compute the gradient, and then update the density variable: 
\begin{lstlisting}
// create a oc optimizer
OCOptimizer oc(0.001, config.designStep, config.dampRatio);
// convergence criteria
ConvergeChecker criteria(config.finthres);
// main loop of optimization
for (int iter = 0; iter < config.max_iter; iter++) {
float val = objective.eval();
// compute derivative
objective.backward(1);
// check convergence
if (criteria.is_converge(iter, val)) { printf("converged\n"); break; }
// make sensitivity symmetry
symmetrizeField(rho.diff(), config.sym);
// update density
oc.update(rho.diff(), rho.value(), config.volRatio);
// make density symmetry
symmetrizeField(rho.value(), config.sym);
}
\end{lstlisting}
	where 
	\lstinline{ConvergeChecker} is a class to check convergence,  \lstinline{symmetrizeField} is a function to symmetrize a 3D tensor according to a given symmetry type.
	In this routine, we do not filter the sensitivity since the density is already filtered.
	\item Output the optimized density field and elastic matrix:
	\begin{lstlisting}
rho.value().toVdb(getPath("rhoFile"));
Ch.writeToTxt(getPath("ChFile"));
	\end{lstlisting}
	where \lstinline{getPath} is a function to prefix the output directory to a given string.
	The member function \lstinline{toVdb} writes the data of \lstinline{TensorVar} to a OpenVDB file.
\end{itemize}

Users only need to define \lstinline{config}, the material interpolation method, and the objective function before running the code to solve the IHPs.
%
The outputs contain
a microstructure visualization file ( \lstinline{*.vbd}), an elastic tensor matrix (\lstinline{*.txt}),
Users can use Rhino to visualize 
\lstinline{*.vdb} files.


\begin{figure*}[!t]
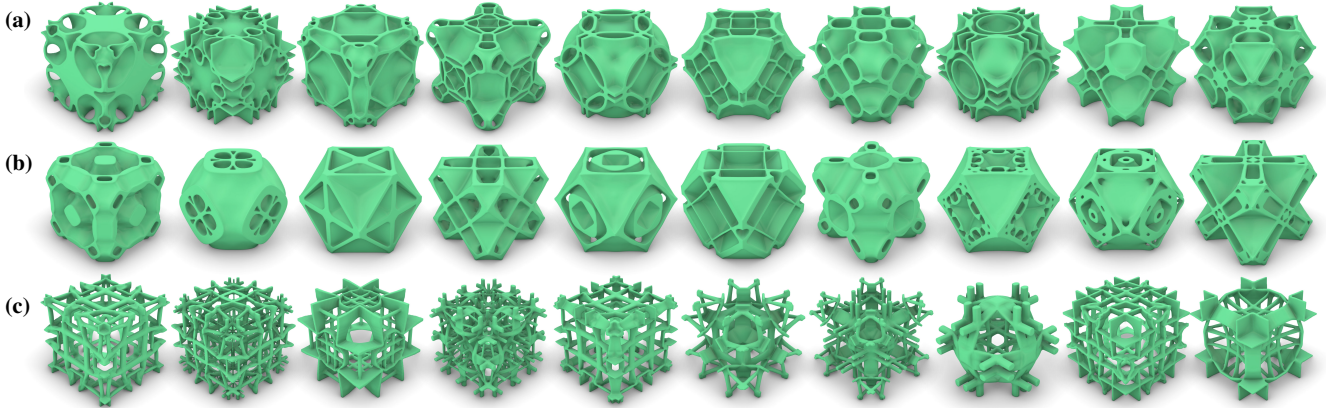

	\centering
	\begin{overpic}[width=0.99\linewidth]{rev-gallery-cc}
		{
			\put(-1.5,8){\textbf{(c)}}
			\put(-1.5,19){\textbf{(b)}}
			\put(-1.5,30){\textbf{(a)}}
		}
	\end{overpic}
	\vspace{-3.5mm}
	\caption{
		Gallery of our optimized microstructures by maximizing bulk modulus (a), maximizing shear modulus (b) and minimizing negative Poisson's ratios (c) with the resolution $128\times 128\times 128$.
	}
	\label{fig:gallery}
\end{figure*}

\subsection{Extensions}
Our framework uses the automatic differentiation (AD) technique~\citep{bookAutoDiff} to make it easy to extend our program to optimize various objective and material interpolation methods.
Users can modify the code according to their needs by changing expressions with different objectives or constraints and material interpolation methods without repeating the tedious calculation.


\paragraph{Different objectives.}
For other objective function, such as shear modulus, its expression can be defined accordingly:
\begin{equation}
	\label{eq:shear-obj}
	f(C^H)=-\frac{1}{3}\left(
	C^H_{33}+C^H_{44}+C^H_{55}
	\right)
\end{equation}
We change nothing than the objective from the code of bulk modulus optimization by calling
\begin{lstlisting}
auto objective = - (Ch(3, 3) + Ch(4, 4) + Ch(5, 5)) / 3.f;
\end{lstlisting}
To design negative Poisson's ratio materials, \cite{xia2015design} propose a relaxed form of objective function for 2D problems.
Accordingly, we can define a similar objective to design negative Poisson's ratio materials in 3D:
\begin{equation}\label{eq:poisson-relx}
	f(C^H)= C^H_{01}+C^H_{02}+C^H_{12} -\beta^l\left(	C^H_{00}+C^H_{11}+C^H_{22}	\right),
\end{equation}
where $\beta\in\left(0,1\right)$ is a user-specified constant and the exponential $l$ is the iteration number. The code is: 
\begin{lstlisting}
auto objective = Ch(0, 1) + Ch(0, 2) + Ch(1, 2) - powf(beta, iter) * (Ch(0, 0)+Ch(1, 1)+Ch(2, 2));
\end{lstlisting}
where \lstinline{beta} is a constant 
in $(0,1)$ and \lstinline{iter} is the iteration number in the main loop of optimization.
We also support common mathematical functions, such as exponential and logarithm functions, to define the expression. 
Several works~\cite{radman2013topological,xia2015design} find that the negative value of Poisson's ratio can reach -1 when the shear modulus is much larger than its bulk modulus. 
Accordingly, we can optimize the following objective function to obtain the negative Poisson's ratio:
\begin{equation}\small
	\label{eq:poisson-new}
	\begin{aligned}
		f(C^H)= & \text{log}(1+\eta(C^H_{01}+C^H_{12}+C^H_{20})/(C^H_{00}+C^H_{11}+C^H_{22}))\\
		& +\tau \left(
		C^H_{00}+C^H_{11}+C^H_{22}
		\right)^\gamma,
	\end{aligned}
\end{equation}
where $\eta, \tau, \gamma$ are three parameters.
In our experiments, \tb{we set $\eta = 0.6, \tau = -1/E_0^\gamma$, and $\gamma = 0.5$, where $E_0$ is the Young's modulus of solids.}
We discuss the difference between \eqref{eq:poisson-relx} and~\eqref{eq:poisson-new} in Section~\ref{sec:diff-objs}.


\paragraph{Different material interpolation methods.}
For the routine of Section~\ref{sec:compile-invoke}, we support other convolution kernels (e.g., linear convolution kernels) for density filtering.
Our program is extensible, and users can define their own convolution kernel.
A more direct material interpolation is defined as:
\begin{lstlisting}
auto rhop = rho.pow(3);
\end{lstlisting}
where we only penalize the density variable by the power of $3$ without filtering.
Accordingly, we should filter the sensitivity before updating the density by \lstinline{OCOptimizer}:
\begin{lstlisting}
oc.filterSens(rho.diff(), rho.value(), config.filterRadius);
\end{lstlisting}

\section{Experiments and Applications}

For the optimization parameters, \tb{the material penalization factor} is 3, the filter radius is 2, the maximum iteration number is 300, the iterative step size of density is 0.05, and the damping factor of the OC method is 0.5.
The optimization is stopped when the relative change of the objective function is less than $0.0005$ for three consecutive iterations.
The cubic domain is discretized with 8-node brick elements. The mechanical properties of solids are Young's modulus \tb{$E_0=1$ and Poisson's ratio $\nu = 0.3$}.

We optimize three different objectives: bulk modulus~\eqref{eq:bulk-obj}, shear modulus~\eqref{eq:shear-obj}, and negative Poisson's ratio~\eqref{eq:poisson-relx} (Fig.~\ref{fig:gallery}).
Table~\ref{tab:statis} summarizes the numerical statistics of all examples.
Our solver consumes \tb{less than 40 seconds} for each iteration with a peak GPU memory of \tb{9 GB} for high-resolution examples with $512^3 \approx 134.2$ million elements. \tb{All experiments are executed on a desktop PC with a 3.6 GHz Intel Core i9-9900K, 32GB of memory, and an NVIDIA GTX 1080Ti graphics card with 11 GB graphics card RAM size.}

\begin{figure}[H]
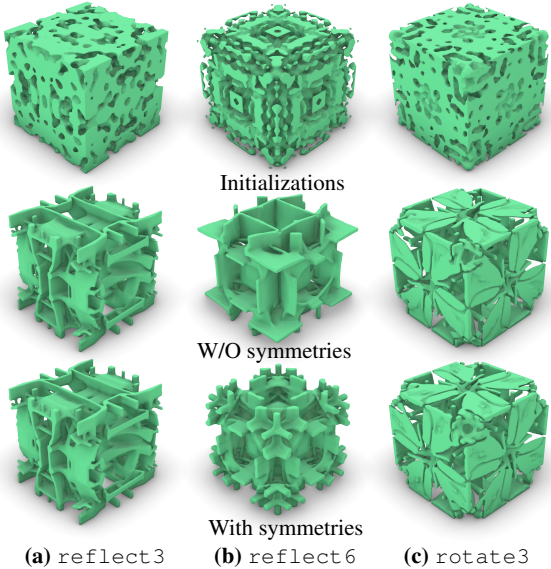

	\centering
	\begin{overpic}[width=0.9\linewidth]{rev-sym-1-cc}
		{
            \put(39,62.5){Initializations}
			\put(35,33){W/O symmetries}
			\put(37,2){With symmetries}
			\put(5,-3){\textbf{(a)} \lstinline{reflect3}}
			\put(38,-3){\textbf{(b)} \lstinline{reflect6}}
			\put(71,-3){\textbf{(c)} \lstinline{rotate3}}
		}
	\end{overpic}
	\vspace{2mm}
	\caption{
		\tb{
			An ablation study of the symmetry operation targeting at negative Poisson's ratio~\eqref{eq:poisson-relx}.
			Upper row: symmetric initializations.
Middle row: without symmetries.
			Bottom row: with symmetries. 
			The resolution is $128\times 128\times 128$ and the volume fraction is 20\%.
		}
	}
	\label{fig:sym}
\end{figure}
\begin{figure*}[!t]
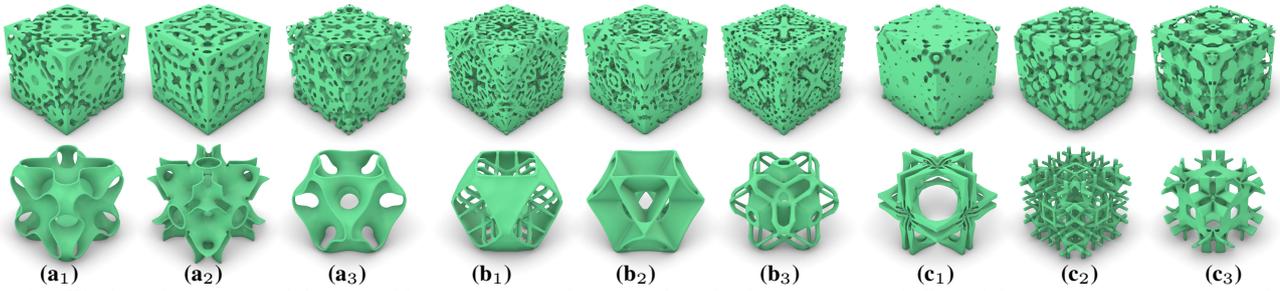

	\centering
	\begin{overpic}[width=0.99\linewidth]{rev-diffinit-cc}
		{
			\put(4,-0){\textbf{(a$_1$)}}
			\put(15,-0){\textbf{(a$_2$)}}
			\put(26,-0){\textbf{(a$_3$)}}
			\put(37,-0){\textbf{(b$_1$)}}
			\put(48,-0){\textbf{(b$_2$)}}
			\put(59,-0){\textbf{(b$_3$)}}
			\put(71,-0){\textbf{(c$_1$)}}
			\put(82,-0){\textbf{(c$_2$)}}
			\put(93,-0){\textbf{(c$_3$)}}
		}
	\end{overpic}
	\vspace{-1mm}
	\caption{
		We show different initial density fields (upper row) and optimized results (bottom row) for bulk modulus maximization (a$_1$)-(a$_3$), shear modulus maximization (b$_1$)-(b$_3$), and negative Poisson's ratio materials (c$_1$)-(c$_3$) with $128\times 128\times 128$ elements under the volume fraction \tb{10\%}.
	}
	\label{fig:diffinit}
\end{figure*}

\subsection{Symmetry}
Symmetry is essential for designing isotropic material.
We have predefined three symmetry types:
\begin{itemize}
\item \lstinline{reflect3}: the reflection symmetry on three planes $\left\{x=0.5,y=0.5,z=0.5\right\}$ of the cube domain;
	\item \lstinline{reflect6}: the reflection symmetry on six planes  $\{x=0.5, y=0.5,z=0.5,x+y=0,y+z=0,z+x=0\}$ of the cube domain;
	\item \lstinline{rotate3}:
	rotation symmetry means that the structure is invariant under the rotation of $90^\circ$ around the x, y, z axes that pass through the cube domain's center, as same under their compositions.
\end{itemize}
Each symmetry splits the cube into many orbits.
To enforce the symmetry, we set the density variables on the same orbit to their average.
Fig.~\ref{fig:sym} shows different symmetry results.
\tb{When the initial value is symmetric, the optimization naturally ensures the symmetry even without reflect3 or rotate3 operation.
However, the operation with or without reflect6 shows the most distinct structural difference. Specifically, the result without using reflect6 only possesses reflect3.
We conjecture that since more symmetry restrictions exist in reflect6 than reflect3 or rotate3, the numerical and machine errors are enlarged with the optimization, thereby causing the symmetric constraint to be violated. Hence, we add symmetry operations to generate symmetric structures.}
In the experiments, we use the \lstinline{reflect6} symmetry by default.

\subsection{Density initializations}
The optimization problem~\eqref{eq:opt-prob} admits a trivial solution, where all the density variables are the predefined volume ratio.
Besides, it~\eqref{eq:opt-prob} has numerous local minima.
The initial density field greatly influences which local minimum it converges to.  
Thus, it is necessary to construct various initial density fields to find desired microstructures.
Previous work usually constructs initialization artificially and seldom discusses other ways of initialization, 
\tb{while we propose to increase the initialization diversity. Different initial density fields and their corresponding optimizated structures are shown in Fig.~\ref{fig:diffinit}. }

We use trigonometric functions to cover various initial density fields.
We first try the following basis functions:
\begin{equation*}
	\begin{aligned}
		T_n=\{
		\cos{2\pi k \bar{x}_i},\sin{2\pi k \bar{x}_i} : 0<k\le n,  i=0,1,2,\\
		\bar{\mathbf{x}}=\mathbf{R}_q\left(\mathbf{x}-\mathbf{b}\right),\mathbf{b}=\left(0.5,0.5,0.5\right)^\top
		\},
	\end{aligned}
\end{equation*}
where the integer $n$ determines the size of the initialization space, 
$\mathbf{x}\in \mathbb{R}^3$ is the coordinate of the element's center, $\mathbf{R}_q\in\mathbb{R}^{3\times3}$ is a rotation matrix determined by a normalized quaternion $q$ with 4 random entries.
Then, to exploit more initializations, we extend $T_n$ as: $Q_n=T_n\cup\left\{p_1p_2 : p_1,p_2\in T_n\right\}$, where
the products of any two items in $T_n$ are incorporated.
$Q_n$ of each element is different.

To initialize a density field, we first generate a set of random numbers in $[-1,1]$ as weights, whose number is the number of the basis functions in $Q_n$.
Then, for each element, we use the obtained weights to weight the basis functions in $Q_n$ and then sum them.
Finally, we project the sum into $[\rho_{\min},1]$ via a rescaled Sigmoid function $S(y)=\rho_{\min} + \widehat{V}/\left({1+e^{-k(y-\mu)}}\right)$, where $k=15$, $\widehat{V} = \min(1.5V,1)$, and $\mu$ is determined by the binary search such that the volume constraint is satisfied after the projection.
This projection aims to produce a valid density distribution, i.e., the constraint $\rho_{min}\leq \rho_{e} \leq 1$ is satisfied for each element, and make the initialization far from the trivial solution.
In Fig.~\ref{fig:diffinit}, we use different initial density fields for optimization.
Different initial density fields lead to different results, which are different local optimal solutions.



\def \g{\cellcolor[HTML]{EFEFEF}}
\linespread{1.08}

\begin{table*}[t]
	\centering
	\caption{\tb{Performance statistics of different precision representations for bulk modulus maximization with $128\times128\times128$ elements under the volume fraction 20\%.  We report the minimum attainable relative residual $r^{\text{min}}_{\text{rel}}$ on the final optimized density field, the memory usage excluding unified memory during optimization (Mem. [MB]), the average time for one iteration (Time/Iter [s]), the whole time cost (Total [s]), and the final bulk modulus (Objective).
			The term ``Sens'' indicates the used memory for sensitivity evaluation.		
			We label the unavailable data with ``-''.}}
	\vspace{-1mm}
	\begin{tabular}{m{1.3cm}|m{1.8cm}|c|c|c|c|c|c|m{1.6cm}|m{1cm}|m{1cm}}
		\hline
		\multirow{2}{1.3cm}{Precision} &\multirow{2}{1.8cm}{ $r_{rel}^{\text{min}}$}  & \multicolumn{6}{c|}{Mem. [MB]} & \multirow{2}{1.6cm}{Time/Iter [s]} & \multirow{2}{1.3cm}{Time [s]}& \multirow{2}{1cm}{Objective}\\
		\cline{3-8}
		& &Density & Stencil & Nodal Vector & Flag & Sensitivity & Total & &\\
		\hline
		FP16      & $1.22\times10^{-2}$ &  8 & 163 & 44 & 8 & 39 & 262 & -  &  - &  - \\ \hline
		FP32      & $2.36\times10^{-6}$ &  8 & 327 & 89 & 8 & 77&  509 & 0.75 & 57 & 0.0678 \\ \hline
		FP64      & $8.01\times10^{-15}$&  8 & 654 & 178 & 8 & 154 & 1002 & 2.05  & 202  & 0.0685 \\ \hline
		FP32/FP16 & $2.13\times{10^{-6}}$&  8 & 163 & 89 & 8 & 0 & 268 & 0.68  & 59  & 0.0684\\ \hline
		FP64/FP32 & $ 8.29\times{10^{-15}}$ &  8 & 327 & 178 & 8 & 0 & 521 & 1.14  & 107  & 0.0685\\ \hline
	\end{tabular}
	\label{tab:pre}
\end{table*}

\subsection{\tb{Mixed-precision scheme}}

\tb{To demonstrate the effectiveness of the proposed mixed-precision approach, we test various precision representations under the same configuration (optimization parameters, resolutions, and desktop PC). The statistics are shown in Table~\ref{tab:pre}.}
\tb{A more precise representation of storage yields a smaller residual, albeit at the cost of increased memory consumption and iteration time. The computation time for pure single precision (FP32) is comparable to that of mixed precision (FP32/FP16). However, utilizing mixed precision (FP32/FP16) can lead to a 47\% reduction comparing with pure FP32 in memory consumption. In addition, the relative error of different precisions in the final bulk modulus is less than 1.1\%.
In summary, due to this mixed-precision scheme, we can solve high-resolution examples with $512^3 \approx 134.2$ million finite elements on only a NVIDIA GeForce GTX 1080Ti GPU.}


\subsection{\tb{Comparison with Multi-CPU framework}}

\tb{We implement the multi-CPU framework~\cite{aage2015topology} and conduct the experiments on a cluster with a total of 9 nodes, each equipped with two Interl Xeon E5-2680 v4 28-core CPUs and 128GB memory connected by Intel OPA.
Since we have verified that the relative residual $10^{-2}$ is acceptable for IHPs (see Fig.~\ref{fig:femrelthres}), the relative residual thresholds for both multi-CPU and our frameworks are set as $10^{-2}$.
In Fig.~\ref{fig:cmp-petsc}, the same initialization is adopted for these two frameworks. The final structures and moduli obtained by both frameworks are very similar.
The average time of each iteration for the Multi-CPU framework is around $40.0$ seconds, while our framework achieves a significantly reduced average time cost of $4.4$ seconds.}

\begin{figure}[t]
	\centering
	\begin{overpic}[width=0.99\linewidth]{cmp-petsc-cc}
		{
			\put(5,-1){Intial density field}
			\put(40,-1){\tb{\small $B = 0.1110$}}
			\put(74,-1){\tb{\small $B = 0.1126$}}
		}
	\end{overpic}
	\vspace{0mm}
	\caption{
	\tb{Maximizing bulk noduli using Multi-CPU framework (Middle) and our framework (Right).
	The domain resolution is set to $256\times256\times256$, and the volume fraction is fixed at $0.3$.}
}
	\label{fig:cmp-petsc}
\end{figure}

\begin{figure}[t]
	\centering
	\begin{overpic}[width=1.0\linewidth]{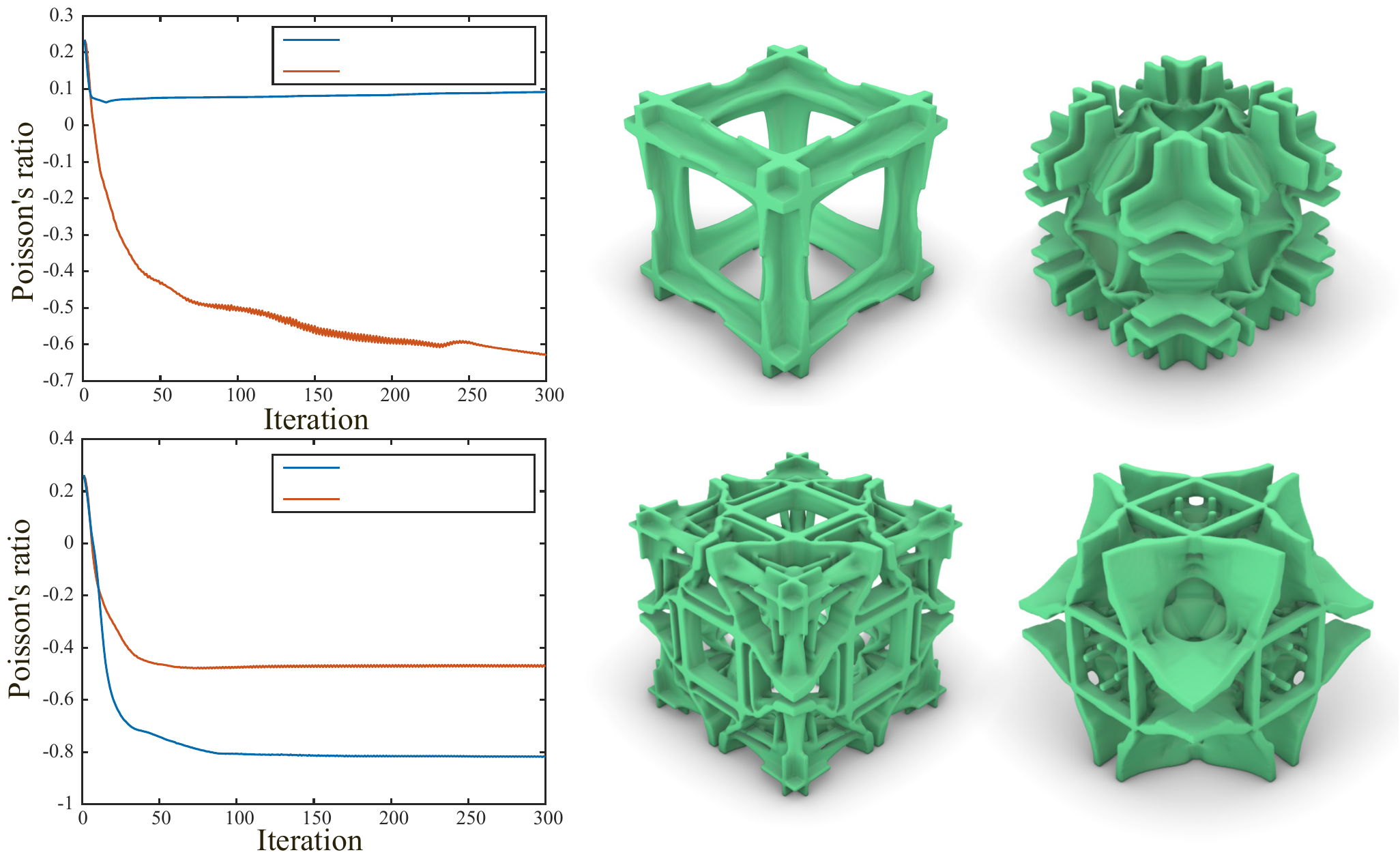}
		{
			\put(24.7,55.4){\tiny Optimizing~\eqref{eq:poisson-relx}}
			\put(24.7,57.4){\tiny Optimizing~\eqref{eq:poisson-new}}
			\put(24.7,24.9){\tiny Optimizing~\eqref{eq:poisson-relx}}
			\put(24.7,27){\tiny Optimizing~\eqref{eq:poisson-new}}
			\put(75,31){\small Optimizing~\eqref{eq:poisson-relx}}
			\put(46,31){\small Optimizing~\eqref{eq:poisson-new}}
			\put(75,0){\small Optimizing~\eqref{eq:poisson-relx}}
			\put(46,0){\small Optimizing~\eqref{eq:poisson-new}}
		}
	\end{overpic}
	\vspace{-4mm}
	\caption{
		Optimizing~\eqref{eq:poisson-relx} and~\eqref{eq:poisson-new} using two different initial density fields.
		The graph plots the Poisson's ratio vs. the number of iterations.
	}
	\label{fig:diff-Poisson-ratio}
\end{figure}
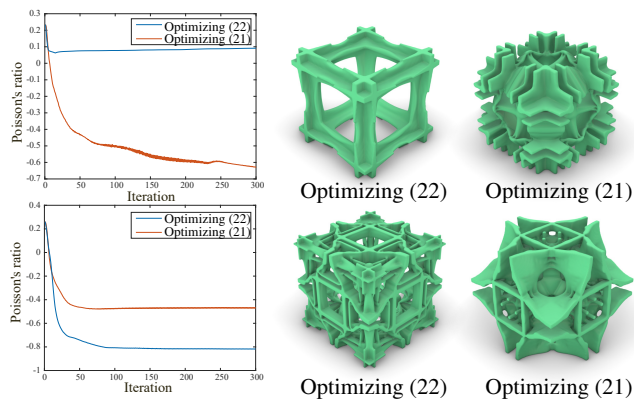

\begin{figure}[t]
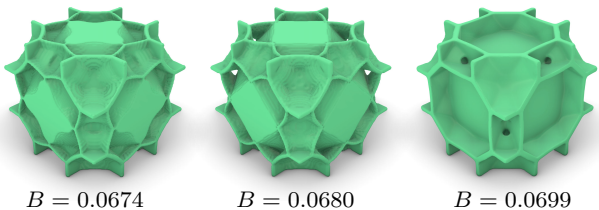

	\centering
	\begin{overpic}[width=1.0\linewidth]{materialinterp-cc}
		{
			\put(7,-1){\small $B = 0.0674$}
			\put(40,-1){\small $B = 0.0680$}
			\put(74,-1){\small $B = 0.0699$}
		}
	\end{overpic}
	\vspace{-4mm}
	\caption{
		\tb{Optimizing the bulk modules $B$ using different material interpolation methods with a same initial density field and a same filtering radius.
		Left: we filter the density field via the \lstinline{Spline4} convolution kernel, i.e., \lstinline{radial_convker_t<float, Spline4>}.
		Middle: we filter the density field by the linear convolution kernel.
		Right: we do not filter the density but filter sensitivity.
		The number of iterations are $282$ (left), $255$ (middle), and $100$ (right), respectively.}
	}
	\label{fig:diff-interp}
\end{figure}

\subsection{Extending our framework}\label{sec:diff-objs}
Users can optimize material properties according to their own goals through our framework.
To verify the scalability of the framework, we optimize \eqref{eq:poisson-relx} and~\eqref{eq:poisson-new} to achieve the negative Poisson's ratio structures using different density fields, as shown in Fig.~\ref{fig:diff-Poisson-ratio}.
From the results, both \tb{objective functions} can lead to negative Poisson's ratios and have their own advantages.
It is an interesting future work to design specific initial density fields so that the objective functions can be optimized to get smaller Poisson's ratios.

In Fig.~\ref{fig:diff-interp}, we test three material interpolation methods for optimizing the bulk modules. 
%
%
The sensitivity filtering of the OC solver is better than density filtering as it has fewer iteration numbers and generates a greater bulk modules.

\begin{figure}[t]
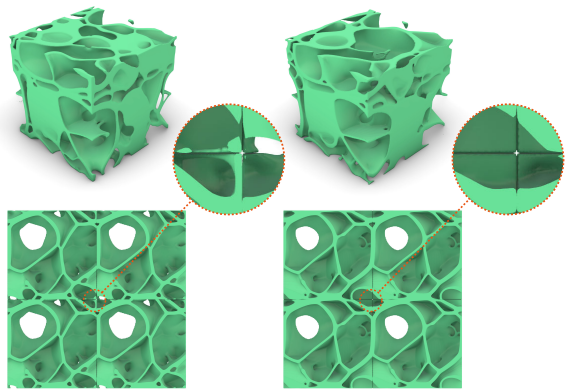

	\centering
	\begin{overpic}[width=0.9\linewidth]{connection-2-cc}
		{
		}
	\end{overpic}
	\vspace{-2mm}
	\caption{
 \tb{Periodic filter kernel.
 We filter sensitivity without (left) and with (right) a periodic filter kernel.}
	}
	\label{fig:diff-filt}
\end{figure}

\tb{
When the symmetry operation is not enforced, the microstructure is not guaranteed to be well connected.
Therefore, we modify the filter kernel to be periodic to improve connectivity (see the zoomed-in views in Fig.~\ref{fig:diff-filt}).
When filtering the sensitivity or density of elements near the boundary, the periodic filter kernel encompasses those elements near the opposite boundary as if multiple unit cells are connected along the boundaries.
}

\subsection{Resolution}
\begin{figure}[t]
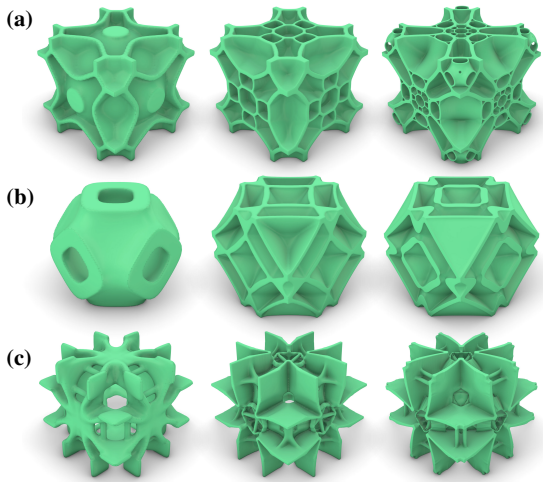

	\centering
	\begin{overpic}[width=0.85\linewidth]{diffreso-cc}
		{
			\put(-3,25){\textbf{(c)}}
			\put(-3,55){\textbf{(b)}}
			\put(-3,88){\textbf{(a)}}
		}
	\end{overpic}
	\vspace{-4mm}
	\caption{
		Various resolutions for bulk modulus maximization (a), shear modulus maximization (b), and negative Poisson's ratio materials (c).
		Left: $64\times64\times64$. Middle: $128\times128\times128$. Right: $256\times256\times256$.
		The volume fraction is $30\%$ for (a \& b) and  $20\%$ for (c).
	}
	\label{fig:diff_reso}
	\vspace{-2mm}
\end{figure}

We test the applications by optimizing \tb{bulk modulus~\eqref{eq:bulk-obj}, shear modulus~\eqref{eq:shear-obj}, and negative Poisson's ratio~\eqref{eq:poisson-relx}} with three different resolutions $64\times 64 \times 64$, $128\times 128 \times 128$, $256\times 256 \times 256$ in Fig.~\ref{fig:diff_reso}.
The increase in computational resolution provided by GPU implementation leads to design improvement.
The respective bulk moduli of the optimized results with different resolutions are \tb{$0.1111$, $0.1124$, and $0.1127$}. The shear moduli of the three structures are \tb{$0.0653$, $0.0721$, and $0.0741$}, respectively. With the increase in resolution, the results show a clear improvement in values and details. A similar conclusion can be obtained for the negative Poisson's ratio microstructures in Fig.~\ref{fig:diff_reso}~(c).

\begin{figure}[t]
	\centering
	\begin{overpic}[width=1\linewidth]{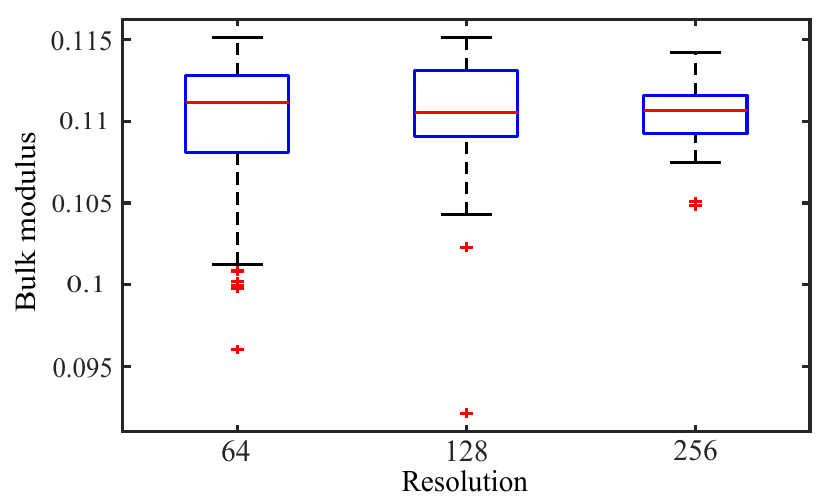}
		{
		}
	\end{overpic}
	\vspace{-2mm}
	\caption{
		Box-plots show the statistics of the optimized bulk moduli for different resolutions ($64\times64\times64$ on the left, $128\times128\times128$ in the middle, and $256\times256\times256$ on the right).
		The volume ratio is \tb{$0.3$}.
	}
	\label{fig:bulk_reso}
	\vspace{-2mm}
\end{figure}

To further validate the effect of resolution on structural properties, we run our optimization on bulk modulus 100 times with different initializations for each resolution, count the resulting bulk modulus, and show the statistics in Fig.~\ref{fig:bulk_reso}. 
The results show that most of the bulk modulus concentrates near the high value when the resolution is high.
\tb{There are also more outliers as the resolution becomes lower. The lower the resolution, the more likely it is to approach the trivial solution.}


\begin{figure}[t]
	\centering
	\begin{overpic}[width=0.85\linewidth]{diffvol-cc}
		{
			\put(-3,25){\textbf{(c)}}
			\put(-3,60){\textbf{(b)}}
			\put(-3,95){\textbf{(a)}}
		}
	\end{overpic}
	\vspace{-5mm}
	\caption{
		Various volume fractions for bulk modulus maximization (a), shear modulus maximization (b), and negative Poisson's ratio materials (c).
		Left: 10\%. Middle: 20\%. Right: 30\%.  
		The resolution is $128^3$.
	}
	\label{fig:diff_vol}
	\vspace{-2mm}
\end{figure}

\begin{figure}[t]
	\centering
	\begin{overpic}[width=0.9\linewidth]{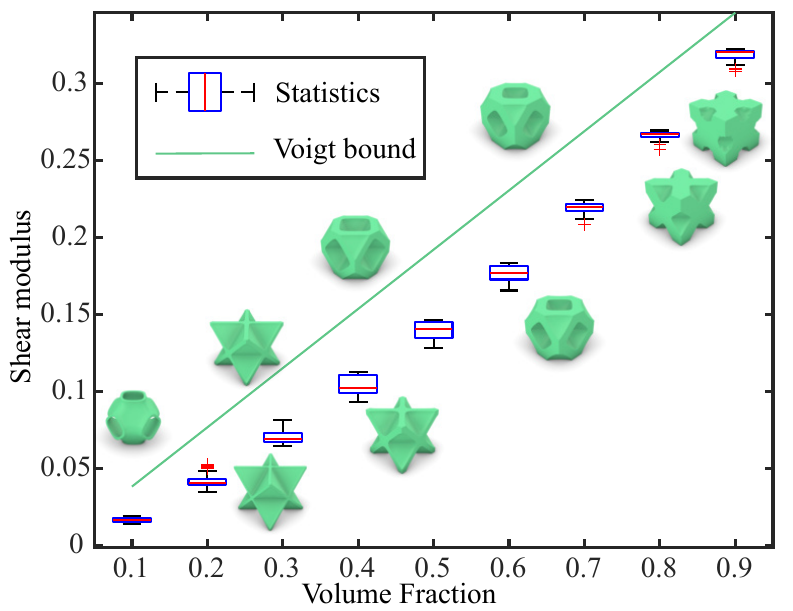}
		{
		}
	\end{overpic}
	\vspace{-2mm}
	\caption{
		\tb{For a $128\times 128 \times 128$  grid, we perform shear modulus maximization 100 times using different initializations for each volume fraction.
		We show statistics of the resulting shear moduli by box-plots.
		The green line shows the Voigt bound.}
	}
	\label{fig:shear_vol}
	\vspace{-2mm}
\end{figure}




\subsection{Volume fraction}
In Fig.~\ref{fig:diff_vol}, three applications are optimized for volume ratios from 10\% to 30\%.
When the volume ratio is 10\%, the results of the bulk modulus and shear modulus optimization are similar to the P surface.
With the increase of volume fraction, the structures become diversified. However, with the increase of volume fraction, the structure changes from rod structure to 
\tb{closed wall structure} for the negative Poisson's ratio structures.
\tb{In addition, we also compare our optimized results with the Voigt bound~\citep{voigt1928lehrbuch}, which provides a theoretical bound of the shear modulus of anisotropic materials under different volume ratios. Increasing the resolution of the microstructure would be considered in the future to obtain the microstructures closer to the upper limit of the theoretical value.} 


\def \g{\cellcolor[HTML]{EFEFEF}}
\begin{table*}[t]
	\centering
	\caption{
		Performance statistics for different objective functions. 
		We report the number of elements (\#Elements), the predefined volume fraction (Vol(\%)), the optimized objective values for three applications (\tb{Objective}),  the peak memory excluding unified memory during optimization (Mem. [MB]), the time for one optimization iteration (Time/Iter [s]), and the whole time cost (Total [s]).
		\tb{The Objective values are calculated after the binarization of the optimized density field.}
}
	\vspace{-1mm}
	\scalebox{0.85}
	{
		\begin{tabular}{c|l|c|c|c|c|c|c|r}
			\hline
			Applications & Examples & Resolution &\#Elements & Vol(\%)  & \tb{Objective}  & Mem. [MB] & Time/Iter [s] & Total [s] \\
			\hline
			\multirow{10}{*}{Bulk modulus} & Fig.~\ref{fig:teaser} (a) & $512^3$ & $1.3 \times 10^8$ & $30$ & $0.1094$ & $8153$ & $38.03$ & $4412$ \\
			& \g Fig.~\ref{fig:diffinit} (a$_1$) & \g $128^3$ & \g $2.1 \times 10^6$ & \g $10$ & \g $0.0239$ & \g $481$ & \g $0.53$ & \g $110$\\
			& Fig.~\ref{fig:diffinit} (a$_2$) & $128^3$ & $2.1 \times 10^6$ & $10$ & $0.0283$ & 481 & $0.51$ & $138$ \\
			& \g Fig.~\ref{fig:diffinit} (a$_3$) & \g $128^3$ & \g $2.1 \times 10^6$ & \g $10$ & \g $0.0301$ & \g $481$ & \g $0.55$ & \g $116$ \\
			& Fig.~\ref{fig:diff-interp} left & $128^3$ & $2.1 \times 10^6$ & $40$ & $0.0674$ & $481$ & $0.55$ & $154$\\
			& \g Fig.~\ref{fig:diff-interp} middle & \g $128^3$ & \g $2.1 \times 10^6$ & \g $40$ & \g $0.0680$ & \g $481$ & \g $0.55$ & \g $141$\\
			& Fig.~\ref{fig:diff-interp} right & $128^3$ & $2.1 \times 10^6$ & $40$ & $0.0699$ & $481$ & $0.64$ & $64$\\
			& \g Fig.~\ref{fig:diff_reso} (a) left & \g $64^3$ & \g $2.6 \times 10^5$ & \g 30 & \g $0.1111$ & \g $197$ & \g $0.15$ & \g $ 12$ \\
			& Fig.~\ref{fig:diff_reso} (a) middle & $128^3$ & $2.1 \times 10^6$ & $30$ & $0.1124$ & $481$ & $0.67$ & $60$\\
			& \g Fig.~\ref{fig:diff_reso} (a) right & \g $256^3$ & \g $1.7 \times 10^7$ & \g $30$ & \g $0.1127$ & \g $2393$ & \g $4.51$ & \g $632$\\
			& Fig.~\ref{fig:diff_vol} (a) left & $128^3$ & $2.1 \times 10^6$ & $10$ & $0.0312$ & $481$ & $0.69$ & $45$ \\
			& \g Fig.~\ref{fig:diff_vol}(a) middle & \g $128^3$ & \g $2.1 \times 10^6$ & \g $20$ & \g $0.0650$ & \g $481$ & \g $0.68$ & \g $50$ \\
			& Fig.~\ref{fig:diff_vol} (a) right & $128^3$ & $2.1 \times 10^6$ & $30$ & $0.1081$ & $481$ & $0.65$ & $64$\\
			\hline
			\multirow{10}{*}{Shear modulus} & \g Fig.~\ref{fig:teaser} (b) & \g $512^3$ & \g $1.3 \times 10^8$ & \g $30$ & \g $0.0684$ & \g $8153$ & \g $36.13$ & \g $5889$ \\
			& Fig.~\ref{fig:diffinit} (b$_1$) & $128^3$ & $2.1 \times 10^6$ & $10$ & $0.0164$ & $481$ & $0.60$ & $45$\\
			& \g Fig.~\ref{fig:diffinit} (b$_2$) & \g $128^3$ & \g $2.1 \times 10^6$ & \g $10$ & \g $0.0181$ & \g $481$ & \g $0.54$ & \g $71$ \\
			& Fig.~\ref{fig:diffinit} (b$_3$) & $128^3$ & $2.1 \times 10^6$ & $10$ & $0.0169$ & $481$ & $0.53$ & $83$\\
			& \g Fig.~\ref{fig:diff_reso} (b) left & \g $64^3$ & \g $2.6 \times 10^5$ & \g $30$ & \g $0.0653$ & \g $197$ & \g $0.15$ & \g $13$ \\
			& Fig.~\ref{fig:diff_reso} (b) middle & $128^3$ & $2.1 \times 10^6$ & $30$ & $0.0721$ & $481$ & $0.60$ & $77$\\
			& \g Fig.~\ref{fig:diff_reso} (b) right & \g $256^3$ & \g $1.7 \times 10^7$ & \g $30$ & \g $0.0741$ & \g $2393$ & \g $4.53$ & \g $671$ \\
			& Fig.~\ref{fig:diff_vol} (b) left & $128^3$ & $2.1 \times 10^6$ & $10$ & $0.0190$ & $481$ & $0.65$ & $62$\\
			& \g Fig.~\ref{fig:diff_vol} (b) middle & \g $128^3$ & \g $2.1 \times 10^6$ & \g 20 & \g $0.0407$ & \g $481$ & \g $0.66$ & \g $76$\\
			& Fig.~\ref{fig:diff_vol} (b) right & $128^3$ & $2.1 \times 10^6$ & $30$ & $0.0804 $ & $481$ & $0.59$ & $84$\\
			\hline
			\multirow{15}{*}{Poisson's ratio} & \g Fig.~\ref{fig:teaser} (c) & \g $512^3$ & \g $1.3 \times 10^8$ & \g $20$ & \g $-0.6644$ & \g $8153$ & \g $34.07$ & \g $8347$ \\
			&  Fig.~\ref{fig:sym} (a) &  $128^3$ &  $2.1 \times 10^6$ &  $20$ & $-0.5768$ &  $481$ &  $0.57$ & $91$\\
			&  \g Fig.~\ref{fig:sym} (b) &  \g $128^3$ &  \g $2.1 \times 10^6$ & \g  20 & \g  $-0.4170$ & \g  $481$ & \g  $0.67$ & \g  $202$ \\
			& Fig.~\ref{fig:sym} (c) & $128^3$ & $2.1 \times 10^6$ & $20$ & $-0.4962$ & $481$ & $0.69$ & $89$ \\
			& \g Fig.~\ref{fig:diffinit} (c$_1$) & \g $128^3$ & \g $2.1 \times 10^6$ & \g $10$ & \g $-0.4399$ & \g $481$ & \g $0.64$ & \g $193$\\
			& Fig.~\ref{fig:diffinit} (c$_2$) &  $128^3$ &  $2.1 \times 10^6$ & $10$ & $-0.5203$ & $481$ & $0.53$ & $64$ \\
			& \g Fig.~\ref{fig:diffinit} (c$_3$) &\g  $128^3$ & \g $2.1 \times 10^6$ & \g $10$ &\g  $-0.4730$ &\g  $481$ &\g  $ 0.50$ & \g $ 150$\\
			&  Fig.~\ref{fig:diff-Poisson-ratio} leftup &  $128^3$ &  $2.1 \times 10^6$ &  $20$ &  $0.0949$ &  $481$ &  $0.49$ &  $147$\\
			& \g Fig.~\ref{fig:diff-Poisson-ratio} leftbottom & \g $128^3$ & \g $2.1 \times 10^6$ & \g $20$ & \g $-0.8130$ & \g $481$ & \g $0.51$ & \g $152$\\
			&  Fig.~\ref{fig:diff-Poisson-ratio} rightup &  $128^3$ &  $2.1 \times 10^6$ &  $20$ &  $-0.6094$ &  $481$ &  $0.58$ &  $175$\\
			& \g Fig.~\ref{fig:diff-Poisson-ratio} rightbottom & \g $128^3$ & \g $2.1 \times 10^6$ & \g $20$ & \g $-0.4363$ & \g $481$ & \g $0.69$ & \g $207$\\
			&  Fig.~\ref{fig:diff_reso} (c) left &  $64^3$ &  $2.6 \times 10^5$ &  $20$ &  $-0.3567$ &  $197$ &  $0.14$ &  $16$\\
			& \g Fig.~\ref{fig:diff_reso} (c) middle & \g $128^3$ & \g $2.1 \times 10^6$ & \g $20$ & \g $-0.5152$ & \g $481$ & \g $0.59$ & \g $126$ \\
			&  Fig.~\ref{fig:diff_reso} (c) right &  $256^3$ &  $1.7 \times 10^7$ &  $20$ &  $-0.6023$ &  $2393$ &  $4.79$ &  $1437$\\
			& \g Fig.~\ref{fig:diff_vol} (c) left & \g $128^3$ & \g $2.1 \times 10^6$ & \g 10 & \g $-0.6301$ & \g $481$ & \g $0.58$ & \g $46$\\
			&  Fig.~\ref{fig:diff_vol} (c) middle &  $128^3$ &  $2.1 \times 10^6$ &  $20$ &  $-0.5587$ &  $481$ &  $0.54$ &  $163$\\
			& \g Fig.~\ref{fig:diff_vol} (c) right & \g $128^3$ & \g $2.1 \times 10^6$ & \g 30 & \g $-0.5470$ & \g $481$ & \g $0.56$ & \g $118$\\
			\hline
		\end{tabular}
	}
	\label{tab:statis}
\end{table*}

\section{Conclusions}\label{sec:conclu}
We have proposed an optimized, easy-to-use, open-source GPU solver for large-scale inverse homogenization problems.
Through a software-level design space exploration, a favorable combination of data structures and algorithms, which makes full use of the parallel computation power of today’s GPUs, is developed to realize a time- and memory-efficient GPU solver.
\tb{Specifically, we use the mixed-precision representation (FP32/FP16) and incorporate padding to handle periodic boundary conditions.
Consequently, this new implementation is carried out on a standard computer with only one GPU operating at the software level.}
%
Topology optimization for achieving high-resolution 3D microstructures becomes computationally tractable with this solver, as demonstrated by our optimized cells with up to $512^3$ (134.2 million) elements.
Our framework is easy-to-use, and the used automatic differentiation technique enables users to design their own objective functions and material interpolation methods.
Code for this paper is publicly available at \url{https://github.com/lavenklau/homo3d}.



\paragraph{Future work and limitations}
Even though our framework is designed to be user-friendly, it needs to modify the source code for specific goals; however, indiscreet modification may produce unexpected compilation or runtime errors.
Heavy dependence on the \tb{templates} makes it harder to track the error.
In future work, we would work on providing a better user interface, e.g., encapsulating the framework as a python module that the user could import.

We support a few material interpolation methods now and will add more \tb{(e.g., RAMP scheme~\cite{stolpe2001alternative}) in future work.}
Besides, we will add support for geometry represented by an implicit function, where the tensor variable becomes the parameters of a set of implicit functions.




\appendix
\section*{Supplementary}
We show the pesudocodes of enfore macro strain Alg.~\ref{alg-force}, Gauss-Seidel relaxation Alg.~\ref{alg-relx}, Assemble numerical stencil for second layer Alg.~\ref{alg-assemble-otf} and evaluate elastic matrix Alg.~\ref{alg-eval-ch} to clearly state the calculation process.


\begin{algorithm}[t]
	\caption{Enforce macro strain}
	\label{alg-force}
	\SetCommentSty{mycommfont}
	\SetKwInput{kwReq}{Require}
	\SetKwInput{kwIn}{Input}
	\SetKwInOut{kwOut}{Output}
	\SetKwFunction{GetNe}{getNeighborElement}
	\SetKwFunction{SetStrain}{setElementMacroStrain}
	\SetKwFunction{GetKE}{loadTemplateMatrix}
	\SetKwFunction{SetLocal}{SetLocal}
	\SetKwFunction{triplerow}{TripleRows}
	\SetKwFunction{BlockIdx}{getBlockIdx}
	\SetKwFunction{threadIdx}{getThreadIdx}
	\SetKwFunction{BlockSize}{getBlockSize}
	\SetKw{kwShare}{Shared}
	\kwIn{Strain id $k$}
	\kwOut{Repsonse force $\mathbf{f}$}
	\kwShare{ $\mathbf{K}_0$} \\
	$\GetKE(\mathbf{K}_0)$\\
	$v= \BlockIdx() \times \BlockSize() + \threadIdx()$ \\
	$\mathbf{f}_v\gets 0 $\\
	\For{$i=0,\cdots,7$} {
		$e \gets \GetNe(vid, i)$\\
		$\mathbf{u}_e\gets \SetStrain(k)$\\
		$\mathbf{f}_v\gets \mathbf{f}_v+ \boldsymbol{\rho}[e]^3\cdot\triplerow(\mathbf{K}_0,7-i)\cdot\mathbf{u}_e$ \tcp*[f]{
			$\triplerow(\mathbf{K}_0,j) $ returns the j-th 3-rows block of matrix $\mathbf{K}_0$.
		}
	}
	$\SetLocal(\mathbf{f},\mathbf{f}_v)$
\end{algorithm}

\def\mtf#1{\text{\textbf{#1}}}
\begin{algorithm}[t]
	\caption{Gauss-Seidel relaxation}
	\label{alg-relx}
	\SetKwInput{kwReq}{Require}
	\SetKwInput{kwIn}{Input}
	\SetKwInOut{kwOut}{Output}
	\SetKwFunction{setBaseId}{setBaseId}
	\SetKwFunction{getIncV}{getIncidentVertex}
	\SetKwFunction{br}{blockReduction}
	\SetKwFunction{GetKE}{loadTemplateMatrix}
	\SetKwFunction{laneId}{getLaneId}
	\SetKwFunction{warpId}{getWarpId}
	\SetKwFunction{gsRelax}{Gauss-SeidelRelax}
	\SetKw{kwShare}{Shared}
	\kwIn{Gauss-Seidel subset id $k$}
	\kwOut{$u$}
	\kwShare{ $\mathbf{K}_0$} \\
	$\GetKE(\mathbf{K}_0)$\\
	$v= \setBaseId(k) + \BlockIdx() \cdot 32 + \laneId()$ \\
	$e\gets \GetNe(v, \warpId())$\\
	\textbf{M}$\gets 0$\\
	\textbf{S}$\gets 0$\\
	\For{$i = 0,1,\cdots,7$}{
		$v^e_i = \getIncV(e, i)$\\
		\eIf{$i\ne 7-\warpId()$}{
			$\mtf{M}\gets\mtf{M}+ \boldsymbol{\rho}[e]^3 \cdot \mathbf{K}^0_{[7-\warpId(),i]}\cdot\mathbf{u}[v^e_i]$
		}{$\mtf{S}\gets\mtf{S}+\boldsymbol{\rho}[e]^3\cdot \mathbf{K}^0_{[i,i]}$}
	}
	$\mtf{M}\gets \br(\mtf{M})$\\
	$\mtf{S}\gets \br(\mtf{S})$\\
	\If{$\warpId()=0$}{
		$\mathbf{u}_v\gets \gsRelax(\mtf{S},\mtf{M},\mathbf{u}[v], \mathbf{f}[v])$\\
		$\mathbf{u}[v]\gets \mathbf{u}_v$
	}
\end{algorithm}

\begin{algorithm}[t]
	\caption{Assemble numerical stencil for the second layer}
	\label{alg-assemble-otf}
	\SetKwInput{kwReq}{Require}
	\SetKwInput{kwIn}{Input}
	\SetCommentSty{mycommfont}
	\SetKwInOut{kwOut}{Output}
	\SetKwFunction{CoveredElements}{CoveredElements}
	\SetKwFunction{isCovered}{isCovered}
	\SetKwFunction{getElementsId}{getElementsId}
	\SetKwFunction{weight}{weight}
	\SetKwFunction{offset}{offset}
	\SetKwFunction{threadIdx}{getThreadIdx}
	\SetKw{kwShare}{Shared}
	\kwOut{Assembled stencil $\mtf{ST}[27]$}
	\kwShare{ $\mathbf{K}_0$} \\
	\GetKE($\mathbf{K}_0$)\\
	$v\gets \BlockIdx() \times \BlockSize() + \threadIdx()$\\
	$v_i\gets 13$\\
	\For{$v_j = 0,1,\cdots,26$}{
		$st\gets 0$\\
		\For{$e_{\text{off}}\mtf{ in } \CoveredElements(v_j)$}{
			\If{$!\isCovered(v_i, e_{\text{off}})$}{continue;}
			$e\gets \getElementsId(v, e_{\text{off}})$\\
			\For{$e\_vi = 0,1,\cdots 8$}{
				$w_i \gets \weight(\offset(v_i,e, e\_vi))$ \tcp*[f]{\eqref{eq:st-assemb-w}}\\
				\For{$e\_vj = 0,1,\cdots 8$}{
					$w_j \gets \weight(\offset(v_j,e, e\_vj))$ \tcp*[f]{\eqref{eq:st-assemb-w}}\\
					$st\gets st+w_i\cdot w_j\cdot \mathbf{K}^0_{[e\_vi,e\_vj]}\cdot \boldsymbol{\rho}[e]^3$
				}
			}
		}
		$\mtf{ST}[v_j][v]\gets st$
	}
\end{algorithm}

\begin{algorithm}
	\caption{Evaluate elastic matrix}
	\label{alg-eval-ch}
	\SetCommentSty{mycommfont}
	\SetKwInput{kwReq}{Require}
	\SetKwInput{kwIn}{Input}
	\SetKwInOut{kwOut}{Output}
	\SetKwFunction{wr}{warpReduction}
	\SetKw{kwShare}{Shared}
	\kwIn{Solved fluctuation displacement $\mathbf{u}[3][nv]$ in \textbf{half2} format}
	\kwOut{Partial sum \textbf{Cs}[21] of \ref{eq:elasttensor}}
	\kwShare{ $\mathbf{K}_0$} \\
	\kwShare{ $\boldsymbol\chi_e[6]$ }\\
	\kwShare{ $\mathbf{u}_e[6][32]$}\\
	
	$\GetKE(\mathbf{K}_0)$\\
	$\SetStrain(\boldsymbol\chi_e)$\\
	$v\gets \BlockIdx()\cdot 32 + \laneId()$\\
	$e\gets \GetNe(v, 0)$\\
	$k\gets\warpId()$\\
	\If{$k<3$}{
		$\mathbf{u}_e^*\gets \mathbf{u}[k][v]$  \tcp*[f]{$\mathbf{u}_e^*$ has type of \textbf{half2} and stores nodal displacements of the element $e$ for two macro strains.} \\
		$\mathbf{u}_e[k*2][\laneId()]\gets \boldsymbol\chi_e[k*2]-\mathbf{u}_e^*.x$\\
		$\mathbf{u}_e[k*2+1][\laneId()]\gets \boldsymbol\chi_e[k*2+1]-\mathbf{u}_e^*.y$
	}
	$\mtf{ce}[21]\gets 0$\\
	$counter\gets0$\\
	\For{$i = 0,1,\cdots,6$}{
		$s\gets 0$\\
		\For{$j=i,i+1,\cdots,6$}{
			$s\gets s + \triplerow(\mathbf{K}_0,k)\cdot \mathbf{u}_e[j][\laneId()]$
		}
		$\mtf{ce}[counter\!+\!+]\gets s^\top\cdot \triplerow(\mathbf{u}_e[i][\laneId()],k)\cdot \boldsymbol{\rho}[e ]^3$\\
	}
	$\br(\mtf{ce})$ \tcp*[f]{The reduction is performed element-wise on $\mtf{ce}$. Those values from the same laneId are summed into the first warp.}\\
	\If{$k=0$}{
		$\wr(\mtf{ce})$ \tcp*[f]{The warp reduction is also performed element-wise on $\mtf{ce}$, and the warp sums all the values with the same offset in \mtf{ce} into the first thread by shuffle operations.} \\
		\If{$\threadIdx()=0$}{
			\For{$i = 0,1,\cdots,21$}{
				$\mtf{Cs}[i][\BlockIdx()]=\mtf{ce}[i]$
			}
		}
		
	}
\end{algorithm}

Considering the code’s readability and ease of use, the code comprises separate classes. The main classes are listed and explained as follows.
\begin{itemize}
	\item[-] \lstinline{Tensor<T>}: This template class indicates a 3D scalar field of type \lstinline{T} on a regular grid with a resolution of $x \times y \times z$, where the
	integers $x, y, z$ are the parameters passed to its constructor. It is
	typically used to store the density field.
	\item[-] \lstinline{TensorVar<T>}: This template class represents a tensor variable. The tensor variable locates in the lowest level of the computation graph and serves as input to other expressions. \lstinline{TensorVar<T>} has two data members of type \lstinline{Tensor<T>}, representing its value and gradient of the objective function, respectively. They are accessed by calling \lstinline{.value()} and \lstinline{.diff()}.
	\item[-] Tensor expression: This is not a class but a family of various classes produced by different mathematical operations
	on \lstinline{TensorVar<T>}. For example, suppose a is of type
	\lstinline{TensorVar<float>}, then \lstinline{1.2*a.pow(3)} is a tensor expression representing the mathematical expression $1.2a^3$, where a corresponds to our tensor variable \lstinline{a}. The tensor expression supports the same mathematical operations as the tensor variable, e.g.,
	we can construct another tensor expression \lstinline{1.2*a.pow(3)+b}
	from \lstinline{1.2*a.pow(3)} and a tensor variable \lstinline{b.} Tensor expressions are used to define material interpolation methods from density variables.
	\item[-] \lstinline{Homogenization}: This class encapsulates the process of numerical homogenization from the density field and sensitivity analysis from a gradient of the elastic matrix. Our dedicated multigrid
	solver is embedded in this class. While it can be used independently, we suggest using it in combination with the template class
	\lstinline{ElasticMatrix<T,Exp>} for simplicity.
	\item[-] \lstinline{ElasticMatrix<Exp>}: This template class possesses two 2D arrays storing the elastic matrix and the gradient of the objective function with respect to the elastic matrix. The constructor of
	\lstinline{ElasticMatrix<Exp>} requires two parameters: (1) an object
	\lstinline{hom of type Homogenization} and (2) an tensor expression object \lstinline{exp} of type \lstinline{Exp} representing a 3D field for the material interpolation. We provide a handy function named \lstinline{genCH(hom,exp)} to construct a \lstinline{ElasticMatrix<Exp>} object, which avoids handling the type of \lstinline{exp}. The usage of \lstinline{ElasticMatrix<Exp>} is to expose the components of the elastic matrix for defining the objective function. We overload the bracket operator () to expose the (i, j)-th component of $C^H$ by calling \lstinline{Ch(i,j)}, where \lstinline{Ch} is an object of \lstinline{ElasticMatrix<Exp>}.
	\item[-] Scalar expression: This expression is similar to the tensor expression, whereas the variable becomes a scalar rather than a tensor. Scalar expression can be used to define the objective function using the exposed component of the elastic matrix, e.g., \lstinline{(Ch(1,1)+Ch(2,2)-Ch(3,3)).pow(2)} defines the objective $(C_{11}^H+C_{22}^H-C_{33}^H)^2$. Scalar expression will not execute computation until the function \lstinline{.eval()} is called. To compute its gradient, we provide a function named \lstinline{.backward(s)}, where \lstinline{s} is a scalar for scaling the gradient. Note that the gradient will be backpropagated until the lowest level of the computation graph is reached, which is usually an object of \lstinline{TensorVar<T>} containing the density variable.
	\item[-] \lstinline{OCOptimizer}: The class implements the Optimality Criteria method. This class is used to filter the sensitivity and update the density variable. The prototype of its constructor is \lstinline{OCOptimizer(float min_density, float stepLimit, float dampExponent)}, where
	\lstinline{min_density} is the minimum density threshold, \lstinline{stepLimit} and \lstinline{dampExponent} are the maximum allowable step and damping ratio when updating the density variable.
\end{itemize}

\section*{Declarations}
\paragraph{Conflict of interest }
The authors declare that they have no conflict of interest.

\paragraph{Replication of results}
Important details for replication of results have been described in the manuscript. Code for this paper is at \url{https://github.com/lavenklau/homo3d}.

\section*{Acknowledgment}
The authors would like to acknowledge the financial support from the Provincial Natural Science Foundation of Anhui (2208085QA01), the Fundamental Research Funds for the Central Universities (WK0010000075), the National Natural Science Foundation of China (61972368 and 62025207).

\bibliography{mybibfile}

\end{document}